\documentclass[11pt,a4paper]{article}
\usepackage{amsfonts}
\usepackage{bbm}
\usepackage{amsmath,amssymb}
\usepackage{mathrsfs}
\usepackage{hyperref}
\usepackage{graphicx}
\usepackage{float}

\newtheorem{thm}{Theorem}[section]
\newtheorem{cor}[thm]{Corollary}
\newtheorem{lem}[thm]{Lemma}
\newtheorem{prop}[thm]{Proposition}

\newtheorem{rem}[thm]{Remark}

\numberwithin{equation}{section}\allowdisplaybreaks

\def\leq{\leqslant}
\def\ge{\geqslant}

\def\leq{\leqslant}
\def\geq{\geqslant}

\begin{document}

\title{\large\bf  Scaling Limit of Modulation Spaces and Their Applications}

\author{\footnotesize \bf Mitsuru Sugimoto$^{\dag}$   and Baoxiang Wang$^{\ddag,}$
\footnote{B. Wang is the corresponding author.  The project is supported in part by NSFC, grant 11771024} \\
\footnotesize
\it $^\dag$Graduate School of Mathematics,  Nagoya University, Chikusa-ku, Nagoya 464-8062, Japan \\
\footnotesize \it $^\ddag$LMAM, School of Mathematical Sciences, Peking University, Beijing 100871, PR of China \\
\footnotesize
\text{ Emails: sugimoto@math.nagoya-u.ac.jp \rm (M.S.) },  wbx@math.pku.edu.cn {\rm (B.W.)}
}
\maketitle

\thispagestyle{empty}
\begin{abstract}
Modulation spaces $M^s_{p,q}$ were introduced by Feichtinger \cite{Fei83} in 1983.  By resorting to the wavelet basis, B\'{e}nyi and Oh \cite{BeOh20} defined a modified version to Feichtinger's modulation spaces for which the symmetry scalings are emphasized for its possible applications in PDE. By carefully investigating the scaling properties of modulation spaces and their connections with the wavelet basis, we consider the scaling limit of modulation spaces, which contains both Feichtinger's and B\'{e}nyi and Oh's modulation spaces.   As their applications, we will give a local well-posedness and a (small data) global well-posedness results for NLS in some scaling limit of modulation spaces,   which generalize the well posedness results of \cite{BeOk09}  and \cite{WaHud07}, and certain super-critical initial data in $H^s$ or in $L^p$ are involved in these spaces. \\

{\bf Keywords:} Modulation spaces, wavelet basis; scaling limit of modulation spaces;  NLS.\\

{\bf MSC 2010:} 42B35, 42B37, 35Q55, 42C40.
\end{abstract}

\section{Introduction}

Let $\mathscr{S}:= \mathscr{S}({\Bbb R}^d)$ be the Schwartz space  and $\mathscr{S}'$ be its dual space. Feichtinger \cite{Fei83} in 1983 introduced the notion of modulation spaces $M^s_{p,q} $ by using the short time Fourier transform (STFT) which is now a basic tool in time-frequency analysis.  Recall that the STFT of a function $f$ with respect to a window function $g \in
\mathscr{S}$ is defined as (see \cite{Fei83, Groch01})
\begin{align} \label{STFT}
V_g f(x, \xi) = \int_{\mathbb{R}^d} e^{-{\rm i} t \cdot \xi}
\overline{g(t-x)} f(t)dt  = \langle f, \ M_\xi T_x g\rangle,
\end{align}
where $T_x:  g \to g(\cdot-x)$ and $  M_\xi : g  \to e^{\mathrm{i}\cdot  \, \xi} g$ denote the translation and modulation operators, respectively. The STFT is closely related to the wave packet transform of C\'{o}rdoba and Fefferman \cite{CoFe78} and the Wigner transform \cite{Groch01}, which is now  a basic tool in time frequency analysis theory.
For any $1\leq p,q\leq \infty$, $s\in \mathbb{R}$,  we denote
\begin{align} \label{amodspace}
\|f\|_{M^s_{p,q}} =  \left\| \langle\xi\rangle ^{s} \|V_g f(x, \xi)\|_{L^p(\mathbb{R}^d_x)}\right\|_{L^q(\mathbb{R}^d_\xi)}.
\end{align}
 Modulation spaces $M^s_{p,q}$ are defined as the
space of all tempered distributions $f\in \mathscr{S}' $ for
which $\|f\|_{M^s_{p,q}}$ is finite (cf. \cite{Fei83}).  Modulation spaces $M^s_{p,q}$ are also said to be Feichtinger spaces.

It is well-known that translation $T_x$, modulation $M_\xi$ and scaling operator $D_\lambda : f\to f(\lambda \, \cdot) $ ($\lambda>0$) are fundamental operators in Banach function spaces. Feichtinger spaces $M^0_{p,q}$ are  invariant under translation and modulation operators. However,  $M^0_{p,q}$ have variant scalings, different functions in $M^0_{p,q}$ may have different scalings and  the sharp
scaling property was obtained in \cite{SuTo07} (see also \cite{HaWa14} for $\alpha$-modulation spaces).
Recently, a modified version of Feichtinger spaces with invariant scalings were introduced by B\'{e}nyi and Oh \cite{BeOh20}. Applying the wavelet basis $\{\psi_{j,k} :=\psi(2^{-j}\cdot  -k)\}_{j\in \mathbb{Z}, k\in \mathbb{Z}^d}$, they introduced the following function spaces $\mathfrak{M}^{\mathbf{w}}_{p,q,r}$ with $p,q,r\in [1,\infty)$. Denote
\begin{align}
& \|f \|_{M_{p,q}^{[j]}} =   \left\| \|\psi_{j,k} (D) f\|_{L^p(\mathbb{R}^d)} \right\|_{\ell^q_k (\mathbb{Z}^d)},   \label{modulation-j}
\end{align}
where
\begin{align}
 \psi_{j,k} (D) f = \mathscr{F}^{-1} \psi_{j,k}   \mathscr{F} f.  \label{psijk}
\end{align}
By carefully choosing the weight $\mathbf{w}=\{w_j\}_{j\in \mathbb{Z}}$, say
\begin{align}
w_j \lesssim \left\{
 \begin{array}{ll}
  2^{-\varepsilon j}, & j\geq 0,\\
  2^{d(1/p+1/q-1 +\varepsilon)j}, & j<0,
\end{array}
\right.  \label{weight}
\end{align}
one can define the following (cf. \cite{BeOh20})
\begin{align}
& \|f \|_{\mathfrak{M}^{\mathbf{w}}_{p,q,r}} =   \left\| w_j\| f\|_{M_{p,q}^{[j]} } \right\|_{\ell^r_j (\mathbb{Z})}.   \label{modulationsymmetry}
\end{align}
For some special weights $\mathbf{w}$, B\'{e}nyi and Oh \cite{BeOh20} showed that $\mathfrak{M}^{\mathbf{w}}_{p,q,r}$ has invariant scalings $\|f(\lambda \, \cdot) \|_{\mathfrak{M}^{\mathbf{w}}_{p,q,r}} \sim \lambda^{a} \|f\|_{\mathfrak{M}^{\mathbf{w}}_{p,q,r}}$ for all $f\in  \mathfrak{M}^{\mathbf{w}}_{p,q,r}$, where $a \in \mathbb{R}$ is independent of $f\in  \mathfrak{M}^{\mathbf{w}}_{p,q,r}$. However, $\mathfrak{M}^{\mathbf{w}}_{p,q,r}$ cannot cover Feichtinger spaces $M^s_{p,q}$, which can be regarded as a modified version of Feichtinger spaces with invariant scalings.

Since $\mathfrak{M}^{\mathbf{w}}_{p,q,r}$ contains two indices $q,r$ to control the growth of the regularity arising from translations and scalings in frequency spaces, the regularity of $\mathfrak{M}^{\mathbf{w}}_{p,q,r}$ depends on both translation index $k$ and scaling index $2^{j}$. It follows that $k$ and $j$ are not completely independent of each other, which leads to that the dual space of $\mathfrak{M}^{\mathbf{w}}_{p,q,r}$ is not  $\mathfrak{M}^{1/\mathbf{w}}_{p',q',r'}$ (see Section \ref{Generalmodspace}).

 A natural relation between wavelet basis and the scalings of Feichtinger spaces is the following dilation identity
 $$
   \| f\|_{M_{p,q}^{[j]} }   = \lambda^{d/p} \| f(\lambda \ \cdot )\|_{M^0_{p,q}}, \ \ \lambda = 2^{-j}.
 $$
By resorting to this identity and the scaling properties of Feichtinger spaces, we will consider the scaling limit spaces of $M^0_{p,q}$, which contain both Feichtinger's and B\'{e}nyi and Oh's modulation spaces. Also, we will characterize the dual spaces of $\mathfrak{M}^{\mathbf{w}}_{p,q,r}$, by a class of new function spaces which are rougher than  $\mathfrak{M}^{1/\mathbf{w}}_{p',q',r'}$.  Moreover, we will study NLS in those scaling limit spaces of $M^0_{p,q}$ and obtain some local and (small data) global well-posedness results,  which generalize the  results of \cite{BeOk09}  and \cite{WaHud07} and certain $L^p$ or $H^s$ super-critical initial data are involved in our results. \\

Throughout this paper, we will use the following notations. $C\ge 1, \ c\le 1$ will denote constants which can be different at different places, we will use $A\lesssim B$ to denote   $A\leqslant CB$; $A\sim B$ means that $A\lesssim B$ and $B\lesssim A$. We write $a\vee b= \max (a,b)$ and $a\wedge b = \min(a,b)$. $\mathbbm{1}_E$ denotes the characteristic function on $E\subset \mathbb{R}^d$.  Let $\Lambda$ be a set with finite many elements, $\# \Lambda$ stands for the number of the elements contained in $\Lambda.$  We denote  $|x|_\infty=\max_{i=1,...,d} |x_i|$, $|x|=|x_1|+...+|x_d|$ for any $x\in \mathbb{R}^d$. $Q(x_0, \delta)$ denotes the cube in $\mathbb{R}^d$ with center at $x_0$ and side length $2\delta$. For any $1\leq p \leq \infty$, we denote by $p'$ the dual number of $p$, i.e., $1/p+ 1/p'=1$.  $L^p=L^p(\mathbb{R}^d)$ ($\ell^p$) stands for the (sequence) Lebesgue space for which the norm is written as $\|\cdot\|_p$ ($\|\cdot\|_{\ell^p}$). The following two Bernstein's inequalities will be frequently used in this paper (cf. \cite{BL76}).

\begin{prop}\label{thm1.9}
{\rm (Bernstein multiplier estimate)}  Let  $1 \leq r\le \infty, \ L \geq [d/2]+1$ and $\rho\in H^L$.   Then
we have
\begin{align} \label{1.35}
\|\mathscr{F}^{-1}\rho \mathscr{F} f\|_r  \lesssim \|\rho\|^{1-d/2L}_{2} \left(\sum^d_{i=1} \|\partial^L_{x_i}\rho \|_{2} \right)^{d/2L} \|f\|_r.
\end{align}
\end{prop}

\begin{prop} [\rm  Proposition 1.3.2 of \cite{Tr83}] \label{thm1.9a}
Let $1\leq p\leq q\leq \infty$, $b>0$, $\xi_0\in \mathbb{R}^d$. Denote $ L^p_{B(\xi_0,b)} = \{f\in L^p: \, {\rm supp}\, \widehat{f} \subset B(\xi_0, b) \}$. Then there exists $C>0$ such that
\begin{align} \label{1.35a}
\|  f\|_q  \lesssim  Cb^{d(1/p-1/q)}  \|f\|_p
\end{align}
holds for all $f\in L^p_{B(\xi_0,b)}$ and $C$ is independent of $b>0$ and $\xi_0\in \mathbb{R}^d$.
\end{prop}

The paper is organized as follows. In Section \ref{Notion} we introduce two kind of of scaling limit spaces of  $M^0_{p,q}$ and show that they are dual spaces in Section \ref{Duality}. In Section \ref{Generalmodspace} we consider a generalized version of the scaling limit spaces of $M^0_{p,q}$ which contain both Feichtinger's and Benyi and Oh's modulation spaces. The dilation property of the scaling limit spaces of $M^0_{p,q}$ will be considered in Section \ref{Scaling}.  The algebraic structure is of importance in applications to PDE and we will obtain an algebraic property of the scaling limit spaces of $M^0_{p,q}$  in Section \ref{AlgebraP}. As applications to NLS, we will obtain a local well-posedness result and a (small data) global well-posedness result in Sections \ref{LNLS} and \ref{GNLS}, respectively.

\section{Scaling limit spaces of $M^0_{p,q}$} \label{Notion}

First, let us recall the equivalent norm on Feichtinger spaces by using the frequency uniform decompositions. Taking notice of   (cf. \cite{Groch01})
\begin{align} \label{STFTprop}
V_g f(x, \xi) = e^{-i x\xi} V_{\hat{g}} \widehat{f} (\xi, -x) = e^{-i x\xi} (\mathscr{F}^{-1} T_{\xi} \overline{\hat{g}} \widehat{f}) (x),
\end{align}
one can consider the frequency-discrete version of the STFT, so-called frequency uniform decomposition operator.  Let $\psi$ be a smooth cut-off function adapted to the unit cube $[-1/2, 1/2]^d$ and $\psi =0$ outside the cube $[-3/4, 3/4]^d$.
We writet $\psi_k =\psi(\cdot - k)$  and assume that
\begin{align} \label{sigmak}
\sum_{k\in \mathbb{Z}^d} \psi_k (\xi) \equiv 1, \  \ \forall \; \
\xi \in
\mathbb{R}^d.
\end{align}
The frequency uniform decomposition operators are defined in the following way:
\begin{align} \label{FUD}
\Box_k := \mathscr{F}^{-1} \psi_k \mathscr{F}, \quad k\in {\Bbb Z}^d.
\end{align}
Feichtinger spaces $M^s_{p,q}({\Bbb R}^d)$  have the following equivalent norms (cf. \cite{WaHud07}):
\begin{align}
  \|f\|_{M^s_{p,q}} :=  \left\|\{\langle k\rangle^s  \Box_k f\}_{k\in \mathbb{Z}^d}  \right\|_{\ell^q(L^p)}. \label{defmod1}
\end{align}
Denote $f_\lambda := f(\lambda \cdot)$.
It is easy to see that
$$
\|\mathscr{F}^{-1} \psi (\lambda \xi -k) \mathscr{F} f\|_p = \lambda^{d/p} \|\Box_k f_\lambda\|_p, \ \ \lambda >0.
$$
Taking the sequence $\ell^q$ norm over all $k\in \mathbb{Z}^d$, we obtain that
\begin{align} \label{dialambda}
\left(\sum_{k\in \mathbb{Z}^d}  \|\mathscr{F}^{-1} \psi (\lambda \xi -k) \mathscr{F} f\|^q_p \right)^{1/q} = \lambda^{d/p} \| f_\lambda\|_{M^0_{p,q}}, \ \ \lambda>0.
\end{align}
Let us start with the dilation property of $M^0_{p,q}$. We have for any $\lambda\geq 1$ (cf. \cite{SuTo07, HaWa14}),
\begin{align}
\lambda^{0 \wedge d\big(\frac1q-\frac1p\big) \wedge d\big(\frac1p+\frac1q-1\big)}\|f\|_{M_{p,q}^0} \lesssim \lambda^{\frac{d}{p}}\|f_\lambda\|_{M_{p,q}^0}  \lesssim   \lambda^{0 \vee d\big(\frac1q-\frac1p\big) \vee d\big(\frac1p+\frac1q-1\big)} \|f\|_{M_{p,q}^0},
\label{dimo}
\end{align}
and for any $0< \nu \leq 1$,
\begin{align}
\nu^{0 \vee d\big(\frac1q-\frac1p\big) \vee d\big(\frac1p+\frac1q-1\big)}\|f\|_{M_{p,q}^0} \lesssim \nu^{\frac{d}{p}}\|f_\nu\|_{M_{p,q}^0}  \lesssim   \nu^{0 \wedge d\big(\frac1q-\frac1p\big) \wedge d\big(\frac1p+\frac1q-1\big)} \|f\|_{M_{p,q}^0}.
\label{2dimo}
\end{align}

\subsection{Refined case $j\leq 0$}

Denote
$$
\psi_{j,k} (\xi):= \psi (2^{-j}\xi - k),  \ \ \Box_{j,k} = \mathscr{F}^{-1} \psi_{j,k} \mathscr{F},   \ \ j\in \mathbb{Z}, \ k\in \mathbb{Z}^d.
$$
Let us write
\begin{align}
\|f \|_{M_{p,q}^{[j]}} =   \left\| \|\Box_{j,k} f\|_{L^p(\mathbb{R}^d)} \right\|_{\ell^q_k (\mathbb{Z}^d)}.   \label{inmodulation-j}
\end{align}
\eqref{dialambda} and \eqref{inmodulation-j} have implied that
\begin{align}
\|f \|_{M_{p,q}^{[j]}} =    \lambda^{d/p} \| f_\lambda\|_{M^0_{p,q}}, \ \ \lambda= 2^{-j}, \ \ j\in \mathbb{Z}.   \label{inmodulation-ja}
\end{align}
Let $\varrho (p,q): \ [1,\infty]\times [1,\infty] \to \mathbb{R}$ and we will write $\varrho:=\varrho (p,q)$ if there is no confusion. We denote by $\mathfrak{M}_{p,q, \infty}^{\varrho}$ the space  of all tempered distributions $f\in \mathscr{S}'$ for which the following norm is finite:
\begin{align}
\|f \|_{\mathfrak{M}_{p,q, \infty}^{\varrho}} = \sup_{j\leq 0} 2^{j\varrho(p,q)}\|f\|_{M_{p,q}^{[j]}}.   \label{invariantmod1}
\end{align}
In order to $\eqref{invariantmod1}$ makes sense, we will always assume that
\begin{align}
 \varrho(p,q) \geq d\left(\frac{1}{p} + \frac{1}{q} -1\right).   \label{invariantmod2}
\end{align}
Using the dilation property in \eqref{dimo},
we have for $j\leq 0$,
\begin{align} \label{dimojleq0}
 2^{-j \left(0 \wedge d\big(\frac1q-\frac1p\big) \wedge d\big(\frac1p+\frac1q-1\big)\right)}  \|f\|_{M_{p,q}^0} \lesssim \|f\|_{M^{[j]}_{p,q}} \lesssim  2^{-j \left(0 \vee d\big(\frac1q-\frac1p\big) \vee d\big(\frac1p+\frac1q-1\big)\right)}  \|f\|_{M_{p,q}^0}.
\end{align}

\begin{lem} \label{basiccondition}
Let $1\leq p,q \leq \infty$.
Then the following statements are equivalent:

\begin{itemize}

\item[\rm (i)] $\varrho (p,q)$ satisfies $\eqref{invariantmod2}$.

\item[\rm (ii)] $  \mathscr{S} \subset \mathfrak{M}_{p,q, \infty}^{\varrho} $.

\item[\rm (iii)] $\mathfrak{M}_{p,q, \infty}^{\varrho} $ contains a nonzero Schwartz function.
\end{itemize}
\end{lem}
{\it Proof.} We have (cf. \cite{BeOh20})
\begin{align} \label{benyiohineq}
\|f\|_{M^{[j]}_{p,q}} \lesssim 2^{-jd(1/p+1/q-1)} \|\langle\nabla \rangle^s f\|_1, \ \ s>d/p.
\end{align}
Hence, if $\eqref{invariantmod2}$ holds then  $\mathscr{S} \subset \mathfrak{M}_{p,q, \infty}^{\varrho}$, i.e., a continuous embedding.  So, condition $\eqref{invariantmod2}$ implies the result of (ii). It suffices to show (iii) $\Rightarrow$ (i). If not, then we have
$$
 \varrho(p,q) <  d\left(\frac{1}{p} + \frac{1}{q} -1\right).
$$
 There exists an $ f \in \mathfrak{M}^{\mu}_{p,q,\infty} \cap \mathscr{S} \setminus \{0\}$.  So, $\widehat{f} \in \mathscr{S} \setminus \{0\}$. Take $\xi_0 \in \mathbb{R}^d$ satisfying $|\widehat{f}(\xi_0)| = \max _{\xi \in \mathbb{R}^d}|\widehat{f}(\xi_0)|$. We can assume, without loss of generality that $\xi_0 =0$. For any $0< \varepsilon \ll 1$, there exists $\delta>0$ such that
$$
|\widehat{f}(\xi) - \widehat{f}(0) | < \varepsilon, \ \ {\rm if} \ \ \xi \in Q(0,\delta)=\{\xi: |\xi|_{\infty} \leq  \delta \}.
$$
Let $j\in \mathbb{Z}_-$ and
$$
\Lambda_j := \{ k\in \mathbb{Z}^d: \ 2^j(k+[-1,1]^d) \subset Q(0,\delta)\}.
$$
It follows that
\begin{align}
\|f\|_{M^{[j]}_{p,q}} \geq & \left( \sum_{k\in \Lambda_j} \|\Box_{j,k}f\|^q_p \right)^{1/q} \nonumber\\
\geq  & \left( \sum_{k\in \Lambda_j} |\widehat{f}(0)|^q \|\mathscr{F}^{-1}\psi_{j,k}  \|^q_p \right)^{1/q}
 - \left( \sum_{k\in \Lambda_j}   \|\mathscr{F}^{-1}[\psi_{j,k} (\widehat{f}(\xi)-\widehat{f}(0) )] \|^q_p \right)^{1/q}.   \label{nonzero1}
\end{align}
It is easy to see that $\# \Lambda_j \sim \delta^d 2^{-jd}$ and so,
\begin{align}
  \left( \sum_{k\in \Lambda_j} |\widehat{f}(0)|^q \|\mathscr{F}^{-1}\psi_{j,k}  \|^q_p \right)^{1/q} \geq c |\widehat{f}(0)| \delta^{d/q} 2^{-jd(1/p+1/q-1)}. \label{nonzero2}
\end{align}
In order to handle the last term in \eqref{nonzero1} is much smaller than that of the right hand side of \eqref{nonzero2}, we divide the proof into the following two cases.

{\it Case 1.} $p>2$. By Hausdorff-Young's inequality,
$$
\|\mathscr{F}^{-1}[\psi_{j,k} (\widehat{f}(\xi)-\widehat{f}(0) )] \|_p \leq \|\psi_{j,k} (\widehat{f}(\xi)-\widehat{f}(0) )\|_{p'} \leq \varepsilon 2^{jd/p'} \|\psi\|_{p'}.
$$
It follows that
\begin{align}
  \left( \sum_{k\in \Lambda_j}   \|\mathscr{F}^{-1}[\psi_{j,k} (\widehat{f}(\xi)-\widehat{f}(0) )] \|^q_p \right)^{1/q} \leq C  \varepsilon \delta^{d/q} 2^{-jd(1/p+1/q-1)} .
  \label{nonzero3}
\end{align}
{\it Case 2.} $1\leq p<2$. We can assume that
$$
|D^\alpha f(\xi)| \leq C, \ \ \xi \in Q(0,\delta).
$$
By H\"older's inequality and Plancherel's identity,
\begin{align}
   \|\mathscr{F}^{-1}[\psi_{j,k} (\widehat{f}(\xi)-\widehat{f}(0) )] \|_p  \leq   \|\psi_{j,k} (\widehat{f}(\xi)-\widehat{f}(0) ) \|^{2-2/p}_2  \|\mathscr{F}^{-1}[\psi_{j,k} (\widehat{f}(\xi)-\widehat{f}(0) )] \|^{2/p-1}_1.
  \label{nonzero4}
\end{align}
One easily sees that
\begin{align}
    \|\psi_{j,k} (\widehat{f}(\xi)-\widehat{f}(0) ) \|_2    \leq \varepsilon 2^{jd/2}\|\psi\|_2.
  \label{nonzero5}
\end{align}
Using Bernstein's inequality, we have for $L:= [d/2]+1$,
\begin{align}
\|\mathscr{F}^{-1}[\psi_{j,k} (\widehat{f}(\xi)-\widehat{f}(0) )] \|_1  & \leq C  \| \psi_{j,k} (\widehat{f}(\xi)-\widehat{f}(0) )  \|^{1-d/2L}_2 \sum^d_{\ell=1} \| \partial^L_{x_\ell}[\psi_{j,k} (\widehat{f}(\xi)-\widehat{f}(0) ) ] \|^{d/2L}_2 \nonumber\\
 & \leq C  (\varepsilon 2^{jd/2})^{1-d/2L}  \sum^d_{\ell=1} \| \partial^L_{x_\ell}[\psi_{j,k} (\widehat{f}(\xi)-\widehat{f}(0) ) ] \|^{d/2L}_2
  \label{nonzero6}
\end{align}
We can assume that $j\ll 0$ with $2^j \leq \varepsilon$. Using $\partial^L (fg) = \sum^L_{l=1} C^l_L \partial^{L-l}f \partial^l g$, we have
$$
|\partial^L_{x_\ell}[\psi_{j,k} (\widehat{f}(\xi)-\widehat{f}(0) ) ]| \leq C  \varepsilon 2^{-jL} \chi_{2^j(k+[-1,1]^d)} (\xi).
$$
It follows that
\begin{align}
  \| \partial^L_{x_\ell}[\psi_{j,k} (\widehat{f}(\xi)-\widehat{f}(0) ) ] \|_2 \leq C  \varepsilon 2^{-jL + jd/2}.
  \label{nonzero7}
\end{align}
By \eqref{nonzero6} and \eqref{nonzero7}
\begin{align}
\|\mathscr{F}^{-1}[\psi_{j,k} (\widehat{f}(\xi)-\widehat{f}(0) )] \|_1   \leq C  \varepsilon .
  \label{nonzero8}
\end{align}
In view of \eqref{nonzero4} and \eqref{nonzero8},
\begin{align}
   \|\mathscr{F}^{-1}[\psi_{j,k} (\widehat{f}(\xi)-\widehat{f}(0) )] \|_p  \leq  C \varepsilon 2^{jd/p'}.
  \label{nonzero9}
\end{align}
So, we have \eqref{nonzero3} in the case $1\leq p <2$.  Collecting the estimates as in \eqref{nonzero1}--\eqref{nonzero3}, we immediately have for $j\ll 0$ and $\varepsilon \ll |\widehat{f}(0)|$,
\begin{align}
 2^{jd(1/p+1/q-1)} \|f\|_{M^{[j]}_{p,q}} \geq c  |\widehat{f}(0)| \delta^{d/q} .   \label{nonzero10}
\end{align}
Hence, for any $\varrho(p,q) < d(1/p+1/q-1)$, we have
$$
  2^{j\varrho (p,q)} \|f\|_{M^{[j]}_{p,q}} \geq c |\widehat{f}(0)|  \delta^{d/q}  2^{j [\varrho(p,q)-d(1/p+1/q-1)] }   .
$$
Taking $j\to -\infty$, we immediately obtain that $\|f \|_{\mathfrak{M}_{p,q, \infty}^{\varrho}}=\infty$. A contradiction. $\hfill\Box$\\

Recall that in the definition of $\mathfrak{M}_{p,q, r}^{\bf w}$ in \eqref{modulationsymmetry}, the weight ${\bf w}$ in \eqref{weight} does not cover the case $w_j = 2^{d(1/p+1/q-1)j}$ for $j\leq 0$.   However,
$\mathfrak{M}_{p,q, \infty}^{\varrho}$ contains $M^0_{p,q}$, we have
\begin{lem} \label{equivalentnorm}
Let $1\leq p,q \leq \infty$ and
\begin{align}
\varrho(p,q) \geq d\left(\frac{1}{p} + \frac{1}{q} -1\right) \vee d\left(\frac{1}{q} - \frac{1}{p} \right) \vee 0.
\end{align}
Then we have $\mathfrak{M}_{p,q, \infty}^{\varrho} = M^0_{p,q}$  with equivalent norm.
\end{lem}
{\it Proof.}
Using the dilation property in \eqref{dimojleq0}, we have for $j\leq 0$,
$$
\|f\|_{M^{[j]}_{p,q}} \lesssim  2^{-j \varrho(p,q)}  \|f\|_{M_{p,q}^0}.
$$
It follows that  $ \|f\|_{\mathfrak{M}_{p,q, \infty}^{\varrho} } \lesssim    \|f\|_{M_{p,q}^0}.$  Taking notice of $M_{p,q}^0= M^{[0]}_{p,q}$, we see that $ \|f\|_{\mathfrak{M}_{p,q, \infty}^{\varrho} } \geq    \|f\|_{M_{p,q}^0}.$ $\hfill\Box$\\

\begin{rem} \label{remclassification}
\rm In view of Lemmas \ref{basiccondition} and \ref{equivalentnorm},  $\mathfrak{M}_{p,q, \infty}^{\varrho}$ has no sense if $\varrho(p,q)< d(1/p+1/q-1)$; $\mathfrak{M}_{p,q, \infty}^{\varrho}$ coincides with $ {M}_{p,q}^{0}$ if $\varrho(p,q) \geq d(1/p+1/q-1) \vee d(1/q-1/p) \vee 0$. In the cases
\begin{align} \label{realsubspace}
 d\left(\frac1p+\frac1q-1 \right) \leq \varrho(p,q)  \leq  d \left(\frac1q- \frac1p \right) \vee 0,
\end{align}
$\mathfrak{M}_{p,q, \infty}^{\varrho}$ seems nontrivial and it is a new class of function spaces. We now introduce a subspace of $\mathfrak{M}_{p,q, \infty}^{\varrho}$:
\begin{align} \label{subspace}
\mathfrak{M}_{p,q, 0}^{\varrho} =\left\{f\in \mathfrak{M}_{p,q, \infty}^{\varrho}: \lim_{j\to -\infty} 2^{j\varrho(p,q)} \|f\|_{M^{[j]}_{p,q}} =0  \right\}.
\end{align}

\begin{figure}
\centering
\includegraphics[width=2.88in, keepaspectratio]{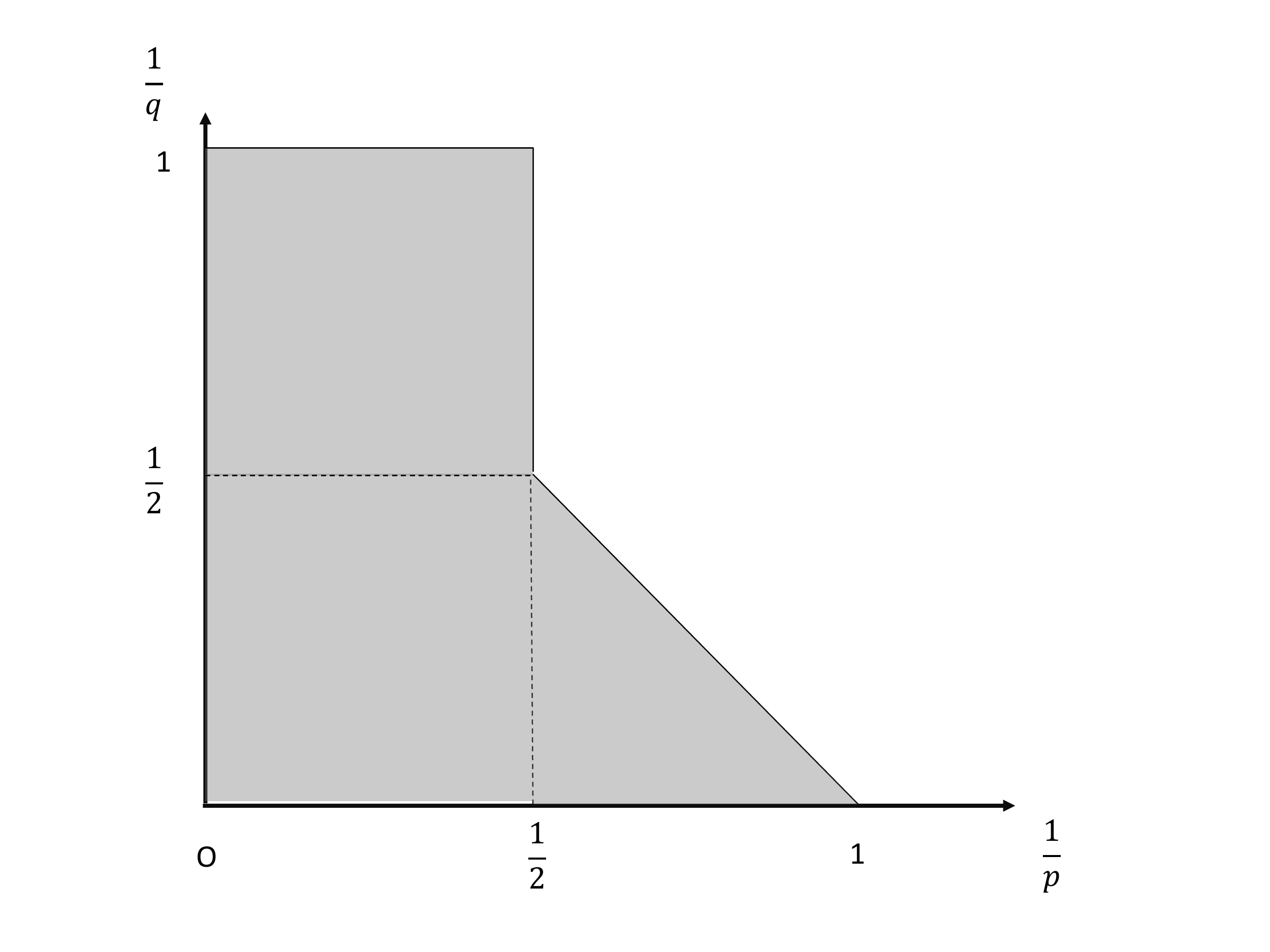}\\
\caption{The region of $(1/p,1/q)$ for which $\varrho(p,q)$ satisfies \eqref{realsubspace}}
\end{figure}
\noindent By \eqref{nonzero10}, any Schwartz function does not belong to $\mathfrak{M}_{p,q, 0}^{\varrho}$ if $\varrho =d(1/p+1/q-1)$. So,  in order to keep that $\mathfrak{M}_{p,q, 0}^{\varrho}$ is a nontrivial Banach space, we need to assume that $\varrho \neq d(1/p+1/q-1)$.
\end{rem}
By the definition of $\mathfrak{M}_{p,q, \infty}^{\varrho}$ and $\mathfrak{M}_{p,q, 0}^{\varrho}$, we have
\begin{prop} \label{Banach-1}
Let $1\leq p,q\leq \infty$ satisfy \eqref{realsubspace}. Then  $\mathfrak{M}_{p,q, \infty}^{\varrho}$ and $\mathfrak{M}_{p,q, 0}^{\varrho}$ are  Banach spaces and $\mathfrak{M}_{p,q, 0}^{\varrho} \subset \mathfrak{M}_{p,q, \infty}^{\varrho} \subset M^0_{p,q}$.
\end{prop}

\begin{rem} \rm Recall that the norm on Wiener amalgam spaces $(L^p, \ell^q)$ can be defined by $\|f\|_{p,q} = \|\|f \psi_{k}\|_{L^p(\mathbb{R}^d)} \|_{\ell^q}$ (see \cite{Wie26}). By considering the scaling limit version of Wiener amalgam spaces, Fofana \cite{Fo88} introduced, so called Fofana spaces $(L^p,\ell^q)^\alpha$ for which the norm can be defined as
$$
\|f\|_{p,q,\alpha}:= {\rm sup }_{j\in \mathbb{Z}} 2^{jd(1/\alpha- 1/q)}\|\|f \psi_{k,j }\|_{L^p(\mathbb{R}^d)} \|_{\ell^q_k}.
$$
We easily see that Wiener amalgam and Feichtinger spaces are a pair of spaces by decomposing the physical and frequency spaces, respectively in a uniform partition manner. Fofana space and $\mathfrak{M}_{p,q, \infty}^{\varrho}$ are a pair of spaces by using the scaling uniform partition to decompose the physical and frequency spaces, respectively.  For each $j$, we see that
$$
 \|f\|_{M^{[j]}_{p,q}} \sim  \|f\|_{M^{0}_{p,q}}, \ \ \ \|\|f \psi_{k,j }\|_{L^p(\mathbb{R}^d)} \|_{\ell^q_k} \sim \|f\|_{p,q}
$$
 and the equivalent constants depend on $j$. So, only the limit behaviour of $j\to  - \infty$ determines that $\mathfrak{M}_{p,q, \infty}^{\varrho}$ is different from Feichtinger space $M^{0}_{p,q}$, from which one can regard that $\mathfrak{M}_{p,q, \infty}^{\varrho}$ is a scaling limit space of $M^0_{p,q}$. Similarly, Fofana space is a scaling limit of Wiener amalgam space.
Recently, the pre-dual spaces of Fofana spaces were studied in \cite{FeFo19,Kp20}.
\end{rem}

 By the inequality of the left hand side of \eqref{dimojleq0}, using an analogous way as in Proposition \ref{equivalentnorm}, we can show that the norm on $M^0_{p,q}$ has another equivalent version:
\begin{align}
\|f\|_{M^0_{p,q}} \sim \inf_{j\leq 0}  2^{j \big(0 \wedge d\big(\frac1q-\frac1p\big) \wedge d\big(\frac1p+\frac1q-1\big)\big)} \|f\|_{M^{[j]}_{p,q}}.  \label{equivinf}
\end{align}
At first glance, one may think that we can define
\begin{align}
\|f\|:  = \inf_{j\leq 0}  2^{j  \mu (p,q)} \|f\|_{M^{[j]}_{p,q}}.  \label{equivinf2}
\end{align}
as another generalization of Feichtinger spaces. Unfortunately, \eqref{equivinf2} is not a norm and so, we are obligated to define
\begin{align}
 \mathscr{M}^{\mu}_{p,q}  = \left\{ f\in \mathscr{S}': \exists \ f_j \in M^{[j]}_{p,q} \ \mbox{such that} \ f= \sum_{j\leq 0} f_j, \ \  \sum_{j\leq 0}    2^{j  \mu } \|f_j\|_{M^{[j]}_{p,q}}<\infty \right\}  \label{equivinf3}
\end{align}
and the norm on $ \mathscr{M}^{\mu}_{p,q}$ is defined as
\begin{align}
\|f\|_{\mathscr{M}^{\mu}_{p,q}} = \inf   \sum_{j\leq 0}    2^{j  \mu } \|f_j\|_{M^{[j]}_{p,q}},   \label{equivinf4}
\end{align}
where the infimum is taken over all of the decompositions of $f= \sum_{j\leq 0} f_j \in \mathscr{M}^{\mu}_{p,q}$.\\

We need to consider the ranges of $\mu:=\mu(p,q)$ into the following three cases.

{\it Case 1.} $\mu(p,q) > d(1/p+1/q-1)$.   We claim that every Schwartz function has zero norm in $\mathscr{M}^{\mu}_{p,q}$. In fact, for any $j_0\in \mathbb{N}$,
we take $f_{j_0}=f \in \mathscr{S}$ and $f_j=0$ for $j\neq j_0$. It follows from \eqref{benyiohineq} that
$$
\|f\|_{\mathscr{M}^{\mu}_{p,q}} \lesssim 2^{j_0(\mu(p,q) - d(1/p+1/q-1))} \|\langle\nabla\rangle^s f\|_{1}.
$$
Taking $j_0\to -\infty$, we have $\|f\|_{\mathscr{M}^{\mu}_{p,q}} =0$.

{\it Case 2.} If $\mu(p,q) \leq 0 \wedge d(1/p+1/q-1) \wedge d(1/q-1/p)$, then we claim that $\mathscr{M}^{\mu}_{p,q} = {M}^{0}_{p,q}$ with equivalent norm.  Indeed, from the inequality of the left hand side of \eqref{dimojleq0} it follows that
$$
2^{\mu(p,q) j}\|f\|_{{M}^{[j]}_{p,q}} \gtrsim  \|f\|_{{M}^{0}_{p,q}}, \ \ j\leq 0.
$$
Hence, we have for and decomposition of $f$, $f= \sum _{j\leq 0} f_j$,
$$
\sum_{j\leq 0} 2^{\mu(p,q) j}\|f_j\|_{{M}^{[j]}_{p,q}} \gtrsim  \sum_{j\leq 0} \|f_j\|_{{M}^{0}_{p,q}} \geq \|f\|_{{M}^{0}_{p,q}},
$$
which implies that $\|f\|_{\mathscr{M}^{\mu}_{p,q}}  \gtrsim  \|f\|_{{M}^{0}_{p,q}}$.  On the other hand, taking $f_0=f$ and $f_j=0$ for any $j\geq 1$, we see that
$$
\sum_{j\leq 0} 2^{\mu(p,q) j}\|f_j\|_{{M}^{[j]}_{p,q}}  \leq \|f_0\|_{{M}^{0}_{p,q}} =  \|f\|_{{M}^{0}_{p,q}}.
$$
Hence, we have $\|f\|_{\mathscr{M}^{\mu}_{p,q}}  \leq  \|f\|_{{M}^{0}_{p,q}}$.

{\it Case 3.} The case $ 0   \wedge d(1/q-1/p) < \mu(p,q)  \leq  d(1/p+1/q-1) $ seems nontrivial. By \eqref{dimojleq0}, we can assume that that
\begin{align} \label{mulimit}
\|f\|_{M^{[j]}_{p,q}} =  2^{-j \mu_j (p,q)} \|f\|_{M^{0}_{p,q}}
\end{align}
for some $  \mu_j(p,q)\in [0\wedge d(1/q-1/p) \wedge d(1/q+1/p -1),  0\vee d(1/q-1/p) \vee d(1/q+1/p -1)]$.  At first glance, one may think that
$\mathscr{M}^{\mu}_{p,q}$ is only a semi-normed space. Assume that
\begin{align} \label{vanishcond}
\mu(p,q) -  \mu_j(p,q)  > \eta >0, \ \ j < -J.
\end{align}
Taking $f_\ell = f$ and $f_j =0$ for any $j \neq \ell $. Then for $\ell \leq -J$,
$$
\sum_{j \leq 0} 2^{\mu(p,q) j} \|f_j\|_{M^{[j]}_{p,q}} \leq    2^{\eta \ell} \|f\|_{M^{0}_{p,q}}   \to 0, \ \ \ell\to -\infty.
$$
It follows that $\|f\|_{\mathscr{M}^{\mu}_{p,q}} =0$.  However, we can show that \eqref{vanishcond} never happens and  have the following

\begin{prop} \label{Banachdense2}
Let $1\leq p,q \leq \infty$ and
\begin{align} \label{dualcondition}
0\wedge d\left(\frac{1}{q} -\frac{1}{p}\right) \leq \mu(p,q)  \leq d\left(\frac{1}{q} +\frac{1}{p}-1\right).
\end{align}
Then $\mathscr{M}^{\mu}_{p,q}$ is a Banach space. Moreover, $\mathscr{S}$ is dense in $\mathscr{M}^{\mu}_{p,q}$ if $p,q<\infty$.
\end{prop}

\begin{figure}
\centering
\includegraphics[width=2.88in, keepaspectratio]{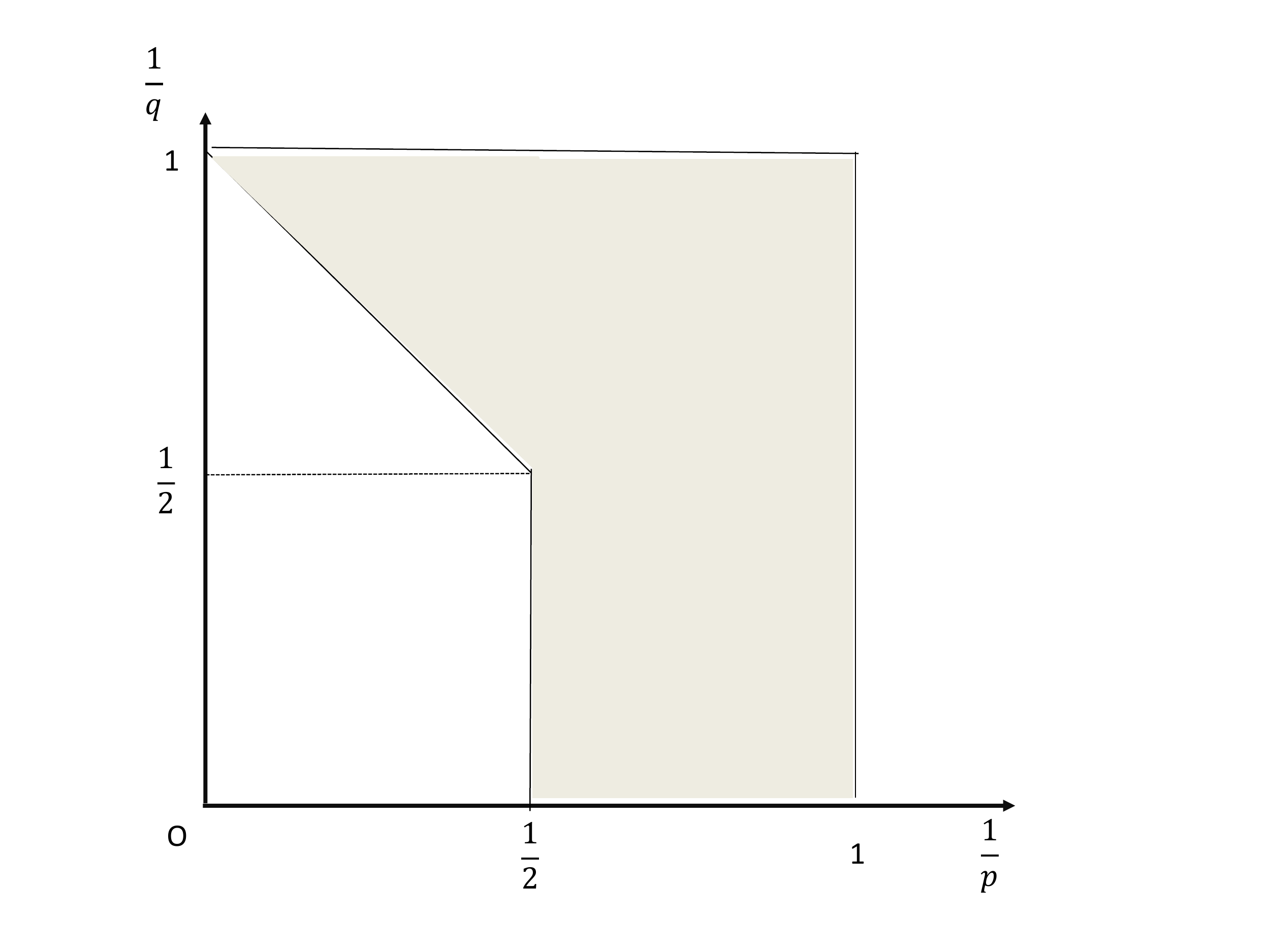}\\
\caption{The grey region of $(1/p,1/q)$ for which $\varrho(p,q)$ satisfies \eqref{dualcondition}}
\end{figure}

{\it Proof.} First, we show that$f\neq 0$ implies $\|f\|_{\mathscr{M}^{\mu}_{p,q}} \neq 0$.
Let $j\in \mathbb{Z}_-$ and
$$
\Lambda_j := \{ k\in \mathbb{Z}^d: \ 2^j(k+[-1/2,1/2]^d)  \cap  Q(0,1) \neq \emptyset\}.
$$
There exists $k_0\in \mathbb{Z}^d$ such that $\Box_{k_0}f \neq 0$ and we can assume that $k_0=0$.  Applying Bernstein's and H\"older's inequalities, we have
\begin{align}
\|\Box_0 f\|_\infty & \leq \sum_{k\in \Lambda_j} \|\Box_{j,k} f\|_\infty \nonumber\\
 & \leq \sum_{k\in \Lambda_j}  2^{jd/p}\|\Box_{j,k} f\|_p \nonumber\\
 &\lesssim  2^{jd(1/p+1/q-1)} \left(\sum_{k\in \Lambda_j} \|\Box_{j,k} f\|^q_p \right)^{1/q} \nonumber\\
 &\lesssim  2^{jd(1/p+1/q-1)} \|f\|_{M^{[j]}_{p,q}}, \label{fnonzero1}
\end{align}
where the hidden constant in the right hand side of \eqref{fnonzero1} is independent of $j\leq 0$ and $f \in  {M}^{[j]}_{p,q}$.  Let $f=\sum_{j\leq 0} f_j$. It follows from \eqref{dualcondition} and \eqref{fnonzero1} that
$$
\|\Box_0 f\|_\infty \leq \sum_{j\leq 0}  \|\Box_0 f_j\|_\infty  \lesssim   \sum_{j\leq 0} 2^{jd(1/p+1/q-1)} \|f_j \|_{M^{[j]}_{p,q}}   \leq  \sum_{j\leq 0} 2^{j\mu} \|f_j \|_{M^{[j]}_{p,q}}.
$$
So, $\|f\|_{\mathscr{M}^{\mu}_{p,q}} \neq 0$.  Hence,  $ \mathscr{M}^{\mu}_{p,q} $ is a normed space. The completeness of  $ \mathscr{M}^{\mu}_{p,q} $ is a direct consequence of $ \mathscr{M}^{\mu}_{p,q} $ as a dual space of $ \mathfrak{M}^{\mu}_{p',q',0} $  in the case  $1<p,q<\infty$ (see Section \ref{Duality}).  For the case $p$ or $q$ equals to $1$, we can give a straightforward proof. Let $\{f_m\}$  be a Cauchy sequence in  $ \mathscr{M}^{\mu}_{p,q} $. One can choose a subsequence still denoted by  $\{f_m\}$ such that
$$
  \|f_{m+1}- f_{m} \|_{ \mathscr{M}^{\mu}_{p,q} } \leq 1/2^m.
$$
It suffices to show that this subsequence has a limit in  $\mathscr{M}^{\mu}_{p,q}$. Let $g_1=f_1$, $g_m = f_{m+1} -f_m$. By the definition of the norm on $\mathscr{M}^{\mu}_{p,q}$, we can find a decomposition of $g_m$, $g_m = \sum_j g_{m,j} \in \mathscr{M}^{\mu}_{p,q} $ such that
$$
 \sum_{j\in \mathbb{Z}_-}  2^{j\mu} \|g_{m,j} \|_{M^{[j]}_{p,q}}   < \|g_m\|_{ \mathscr{M}^{\mu}_{p,q} } + 1/2^m.
$$
It follows that
\begin{align} \label{complet-a}
\sum^\infty_{m=1} \sum_{j\in \mathbb{Z}_-}  2^{j\mu} \|g_{m,j} \|_{M^{[j]}_{p,q}}   < \sum_{m\geq 1}\|g_m\|_{ \mathscr{M}^{\mu}_{p,q} } + 1<2.
\end{align}
 Taking $f_j = \sum_{m\geq 1} g_{m,j}$. By the completeness of $M^{[j]}_{p,q}$, we have $f_j \in M^{[j]}_{p,q}$.  Put $f=\sum_{j\leq 0} f_j = \sum_{j\leq 0} \sum_{m\geq 1} g_{m,j}$. The above inequality implies that $f\in \mathscr{M}^{\mu}_{p,q}$.  Let us observe that
 $$
 f- f_{n+1} = \sum_{m\geq 1} g_m - \sum^n_{m=1} g_m = \sum_{m>n} g_m =\sum_{m>n}\sum_{j\leq 0} g_{m,j},
 $$
\eqref{complet-a} has implied that $\|f-f_{n+1}\|_{\mathscr{M}^{\mu}_{p,q}} \to 0$ as $n\to \infty$.

Now we show the density of $\mathscr{S}$ in $ \mathscr{M}^{\mu}_{p,q} $. For $f=\sum_{j\leq 0} f_j \in \mathscr{M}^{\mu}_{p,q}$, there exists $J\in \mathbb{N}$  satisfying
$$
\sum_{j\leq -J} 2^{j\mu} \|f_j\|_{M^{[j]}_{p,q}} < \varepsilon.
$$
 Let us recall that $M^0_{p,q} = M^{[j]}_{p,q}$ for any $j\leq 0$  and  $\mathscr{S}$ is dense in in $  {M}^{0}_{p,q} $. So, one can find $\varphi_j \in \mathscr{S}$ verifying
 $$
  2^{j\mu} \|f_j - \varphi_j\|_{M^{[j]}_{p,q}} < \varepsilon/J,  \ \  j=-J+1, ..., 0.
 $$
Taking $\varphi_j=0$ for $j\leq -J$ and $\varphi = \sum_{j\leq 0} \varphi_j$, we have
$$
\|f-\varphi\|_{\mathscr{M}^{\mu}_{p,q} }  \leq  \sum_{j\leq 0}   2^{j\mu} \|f_j - \varphi_j\|_{M^{[j]}_{p,q}}  <2\varepsilon.
$$
Since $\varepsilon>0$ is arbitrarily small, we get the density of  $\mathscr{S}$ in $ \mathscr{M}^{\mu}_{p,q} $.
$\hfill\Box$

\subsection{Rough case $j\geq 0$}
Using the dilation property in \eqref{2dimo},
we have for $\nu = 2^{-j}$, $j\geq 0$,
\begin{align} \label{2dimojleq0}
 2^{-j \left(0 \vee d\big(\frac1q-\frac1p\big) \vee d\big(\frac1p+\frac1q-1\big)\right)}  \|f\|_{M_{p,q}^0} \lesssim \|f\|_{M^{[j]}_{p,q}} \lesssim  2^{-j \left(0 \wedge d\big(\frac1q-\frac1p\big) \wedge d\big(\frac1p+\frac1q-1\big)\right)}  \|f\|_{M_{p,q}^0}.
\end{align}
Let  $\varrho:=\varrho (p,q) \in \mathbb{R}$.  Let $\tilde{\mathfrak{M}}_{p,q, \infty}^{\varrho}$ stand for the space of all tempered distribution $f\in \mathscr{S}'$ for which the following norm is finite:
\begin{align}
\|f \|_{\tilde{\mathfrak{M}}_{p,q, \infty}^{\varrho}} = \sup_{j\geq 0} 2^{j\varrho(p,q)}\|f\|_{M_{p,q}^{[j]}}.   \label{2invariantmod1}
\end{align}
\begin{lem} \label{2basiccondition}
Let $1\leq p,q \leq \infty$. Then the following are equivalent:
\begin{itemize}
\item[\rm (i)] $\varrho (p,q)\leq 0$.

\item[\rm (ii)] $  \mathscr{S} \subset \tilde{\mathfrak{M}}_{p,q, \infty}^{\varrho} $.

\item[\rm (iii)] $\tilde{\mathfrak{M}}_{p,q, \infty}^{\varrho} $ contains a nonzero Schwartz function.
\end{itemize}
\end{lem}
{\it Proof.} In view of (cf. \cite{BeOh20})
\begin{align} \label{2benyiohineq}
\|f\|_{M^{[j]}_{p,q}} \lesssim   \|\langle\nabla \rangle^s f\|_p, \ j\geq 0,  \ s>d/q,
\end{align}
we see that $ \mathscr{S} \subset \tilde{\mathfrak{M}}_{p,q, \infty}^{\varrho}$ if $\varrho (p,q)\leq 0$.  So, it sufices to show that ${\rm (iii)} \Rightarrow {\rm (i)} $. Let $f\in \mathscr{S} $. Let $\sigma(\xi)=1$ for $|\xi|_\infty \leq 1$ and ${\rm supp} \sigma \subset [-2,2]^d$, $\sigma_{j} =\sigma (2^{-j}\ \cdot)$. Then we have
$$
\|f\|_{M_{p,q}^{[j]}} \gtrsim  \|\sigma_j(D) f\|_p  \geq \|f\|_p -  \|(I-\sigma_j(D)) f\|_p\geq \frac12 \|f\|_p
$$
if $j$ is sufficiently large. It follows that
$$
\sup_{j\geq 0} 2^{j\varrho }\|f\|_{M_{p,q}^{[j]}} =\infty
$$
for any $\varrho >0$.  $\hfill\Box$

\begin{lem} \label{2equivalentnorm}
Let $1\leq p,q \leq \infty$ and
\begin{align}
\varrho(p,q) \leq d\left(\frac{1}{p} + \frac{1}{q} -1\right) \wedge d\left(\frac{1}{q} - \frac{1}{p} \right) \wedge  0.
\end{align}
Then we have $\tilde{\mathfrak{M}}_{p,q, \infty}^{\varrho} = M^0_{p,q}$  with equivalent norm.
\end{lem}

{\it Proof.}
Using the dilation property in \eqref{2dimojleq0}, we have for $j\geq 0$,
$$
2^{j \varrho(p,q) } \|f\|_{M^{[j]}_{p,q}} \lesssim  \|f\|_{M_{p,q}^0}.
$$
It follows that  $ \|f\|_{\tilde{\mathfrak{M}}_{p,q, \infty}^{\varrho} } \lesssim    \|f\|_{M_{p,q}^0}.$  Obviously,  $ \|f\|_{\tilde{\mathfrak{M}}_{p,q, \infty}^{\varrho} } \geq    \|f\|_{M_{p,q}^0}.$ $\hfill\Box$\\

Hence,  In view of Lemmas \ref{2basiccondition} and \ref{2equivalentnorm}, the nontrivial case is
\begin{align} \label{2realsubspace}
d\left(\frac{1}{p} + \frac{1}{q} -1\right) \wedge d\left(\frac{1}{q} - \frac{1}{p} \right)   \leq \varrho(p,q) \leq 0.
\end{align}

\begin{figure}
\centering
\includegraphics[width=2.88in, keepaspectratio]{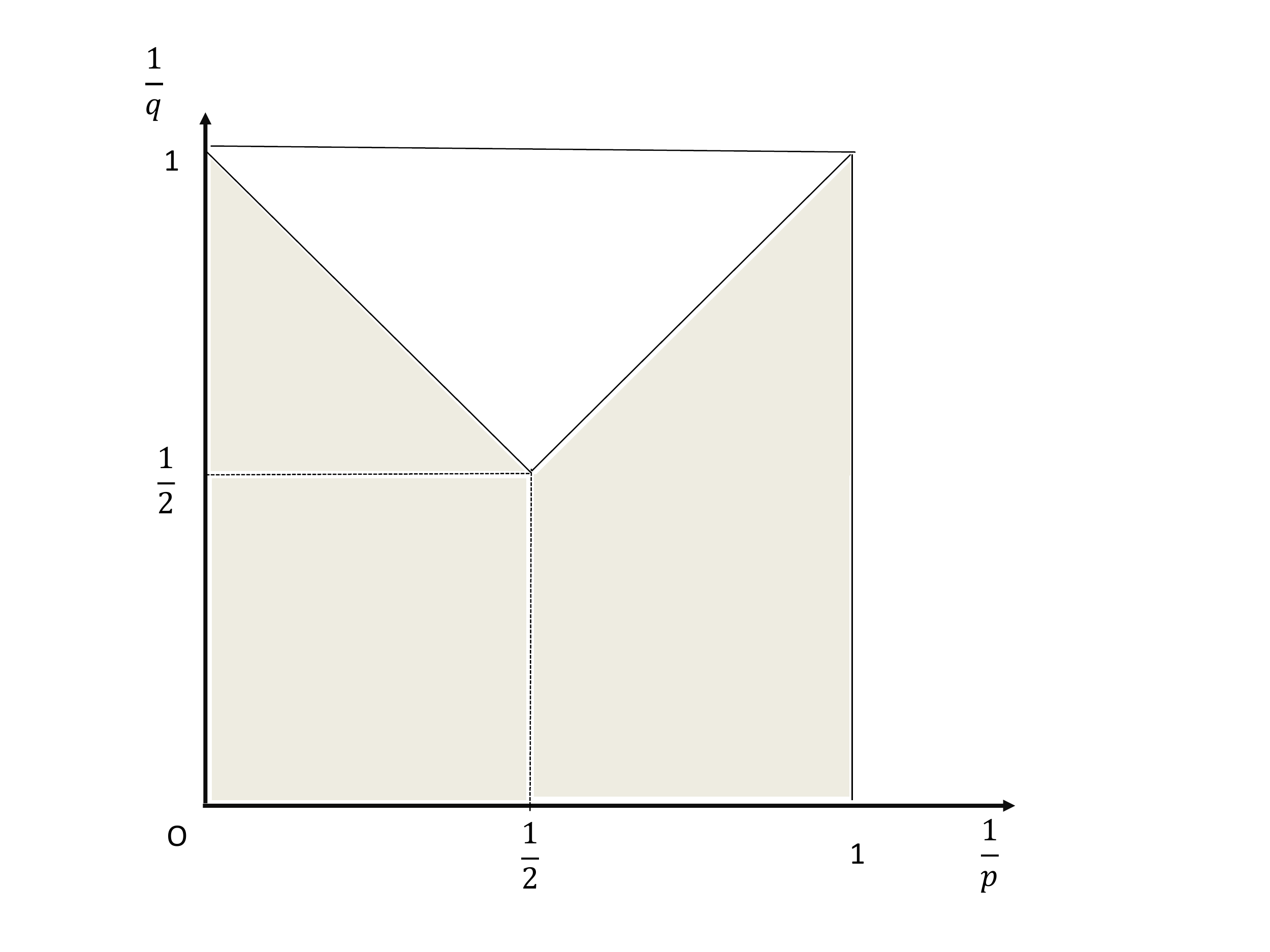}\\
\caption{The grey region of $(1/p,1/q)$ for which $\varrho(p,q)$ satisfies \eqref{2realsubspace}}
\end{figure}

We now introduce a subspace of $\tilde{\mathfrak{M}}_{p,q, \infty}^{\varrho}$:
\begin{align} \label{2subspace}
\tilde{\mathfrak{M}}_{p,q, 0}^{\varrho} =\left\{f\in \tilde{\mathfrak{M}}_{p,q, \infty}^{\varrho}: \lim_{j\to  +\infty} 2^{j\varrho(p,q)} \|f\|_{M^{[j]}_{p,q}} =0  \right\}.
\end{align}
From the proof of Lemma \ref{2basiccondition}, we see that any Schwartz function is not in $\tilde{\mathfrak{M}}_{p,q, 0}^{\varrho}$ if $\varrho=0$. So, in order to keep $\tilde{\mathfrak{M}}_{p,q, 0}^{\varrho}$  is nontrivial, one needs to assume that $\varrho \neq 0$ in $\tilde{\mathfrak{M}}_{p,q, 0}^{\varrho}$.   By the definition of $\tilde{\mathfrak{M}}_{p,q, \infty}^{\varrho}$ and $\tilde{\mathfrak{M}}_{p,q, 0}^{\varrho}$, we have
\begin{prop}
Let $1\leq p,q\leq \infty$ and $\varrho(p,q)$ satisfies \eqref{2realsubspace}. Then  $\tilde{\mathfrak{M}}_{p,q, \infty}^{\varrho}$ and $\tilde{\mathfrak{M}}_{p,q, 0}^{\varrho}$ are  Banach spaces and $\tilde{\mathfrak{M}}_{p,q, 0}^{\varrho} \subset \tilde{\mathfrak{M}}_{p,q, \infty}^{\varrho} \subset M^0_{p,q}$.
\end{prop}

By the first inequality in \eqref{2dimojleq0}, using a similar way as in Lemma \ref{2equivalentnorm}, we can show that the norm on $M^0_{p,q}$ has another equivalent version:
\begin{align}
\|f\|_{M^0_{p,q}} \sim \inf_{j\geq 0}  2^{j \big(0 \vee d\big(\frac1q-\frac1p\big) \vee d\big(\frac1p+\frac1q-1\big)\big)} \|f\|_{M^{[j]}_{p,q}}.  \label{2equivinf}
\end{align}
Since \eqref{2equivinf} is not a norm, we need to define
\begin{align}
 \tilde{\mathscr{M}}^{\mu}_{p,q}  = \left\{ f\in \mathscr{S}': \exists \ f_j \in M^{[j]}_{p,q} \ \mbox{such that} \ f= \sum_{j\geq 0} f_j, \ \  \sum_{j\geq 0}    2^{j  \mu } \|f_j\|_{M^{[j]}_{p,q}}<\infty \right\}  \label{2equivinf3}
\end{align}
and the norm on $ \tilde{\mathscr{M}}^{\mu}_{p,q}$ is defined as
\begin{align}
\|f\|_{\tilde{\mathscr{M}}^{\mu}_{p,q}} = \inf   \sum_{j\geq 0}    2^{j  \mu } \|f_j\|_{M^{[j]}_{p,q}},   \label{2equivinf4}
\end{align}
where the infimum is taken over all of the possible decompositions of $f= \sum_{j\geq 0} f_j \in \tilde{\mathscr{M}}^{\mu}_{p,q}$.\\

{\it Case 1.} $\mu(p,q) < 0$.   We claim that every Schwartz function has zero norm in $\mathscr{M}^{\mu}_{p,q}$. In fact, for any $i \in \mathbb{N}$,
we take $f_{i}=f \in \mathscr{S}$ and $f_j=0$ for $j\neq i$. It follows from \eqref{2benyiohineq} that
$$
\|f\|_{\mathscr{M}^{\mu}_{p,q}} \lesssim 2^{i \mu(p,q)} \|\langle\nabla\rangle^s f\|_{p}, \ \ s>d/q.
$$
Taking $i \to  \infty$, we have $\|f\|_{\mathscr{M}^{\mu}_{p,q}} =0$.

{\it Case 2.} If $\mu(p,q) \geq 0 \vee d(1/p+1/q-1) \vee d(1/q-1/p)$, then we claim that $\tilde{\mathscr{M}}^{\mu}_{p,q} = {M}^{0}_{p,q}$ with equivalent norm.  Indeed, from the inequality of the left hand side of \eqref{2dimojleq0} it follows that
$$
2^{\mu(p,q) j}\|f\|_{{M}^{[j]}_{p,q}} \gtrsim  \|f\|_{{M}^{0}_{p,q}}, \ \ j\geq 0.
$$
Hence, we have for any possible decomposition of $f$, $f= \sum _{j\leq 0} f_j$,
$$
\sum_{j\geq 0} 2^{\mu(p,q) j}\|f_j\|_{{M}^{[j]}_{p,q}} \gtrsim  \sum_{j\geq 0} \|f_j\|_{{M}^{0}_{p,q}} \geq \|f\|_{{M}^{0}_{p,q}},
$$
which implies that $\|f\|_{\tilde{\mathscr{M}}^{\mu}_{p,q}}  \gtrsim  \|f\|_{{M}^{0}_{p,q}}$.  On the other hand, taking $f_0=f$ and $f_j=0$ for any $j\geq 1$, we see that
$$
\sum_{j\geq 0} 2^{\mu(p,q) j}\|f_j\|_{{M}^{[j]}_{p,q}}  \leq \|f_0\|_{{M}^{0}_{p,q}} =  \|f\|_{{M}^{0}_{p,q}}.
$$
Hence, we have $\|f\|_{\tilde{\mathscr{M}}^{\mu}_{p,q}}  \leq  \|f\|_{{M}^{0}_{p,q}}$.

{\it Case 3.} The case $ 0   \leq \mu(p,q)  \leq d(1/q-1/p) \vee  d(1/p+1/q-1) $ seems nontrivial.  We can show that

\begin{prop} \label{2Banachdense2}
Let $1\leq p,q \leq \infty$ and
\begin{align} \label{2dualcondition}
0 \leq \mu(p,q) \leq d\left(\frac{1}{q} -\frac{1}{p}\right) \vee d\left(\frac{1}{q} +\frac{1}{p}-1\right).
\end{align}
Then $\tilde{\mathscr{M}}^{\mu}_{p,q}$ is a Banach space. Moreover, $\mathscr{S}$ is dense in $\tilde{\mathscr{M}}^{\mu}_{p,q}$ if $p,q<\infty$.
\end{prop}

\begin{figure}
\centering
\includegraphics[width=2.88in, keepaspectratio]{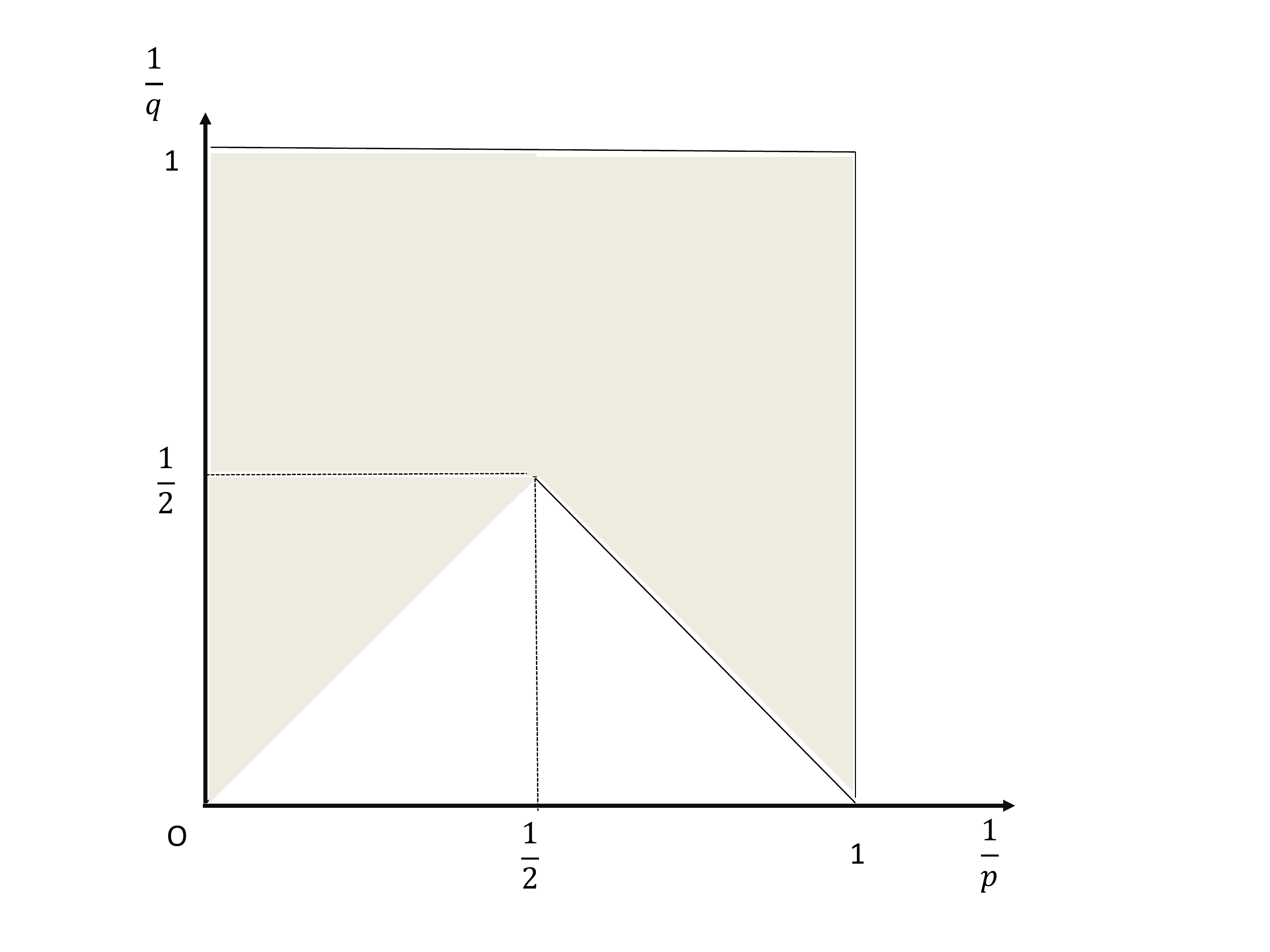}\\
\caption{The grey region of $(1/p,1/q)$ for which $\varrho(p,q)$ satisfies \eqref{2dualcondition}}
\end{figure}

{\it Proof.} First,
There exists $k_0\in \mathbb{Z}^d$ such that $\Box_{k_0}f \neq 0$ and we can assume that $k_0=0$.  We have for any $j\geq 0$,
\begin{align}
\|\Box_0 f\|_\infty  \lesssim \sum_{|\ell|_\infty \lesssim 1} \|\Box_{j,\ell} f\|_\infty.   \label{2fnonzero1}
\end{align}
So, $\|f\|_{\tilde{\mathscr{M}}^{\mu}_{p,q}} \neq 0$.  Hence,  $ \mathscr{M}^{\mu}_{p,q} $ is a normed space. The completeness of  $ \mathscr{M}^{\mu}_{p,q} $ will be shown in the next section.
The density of $\mathscr{S}$ in $ \tilde{\mathscr{M}}^{\mu}_{p,q} $ follows the proof of Proposition \ref{Banachdense2}.
$\hfill\Box$

\subsection{Full case $j\in \mathbb{Z}$}

Let  $\varrho:=\varrho (p,q) \in \mathbb{R}$. We denote by $\dot{\mathfrak{M}}_{p,q, \infty}^{\varrho}$ the space consisting of all tempered distribution $f\in \mathscr{S}'$ for which the following norm is finite:
\begin{align}
\|f \|_{\dot{\mathfrak{M}}_{p,q, \infty}^{\varrho}} = \sup_{j\in \mathbb{Z}} 2^{j\varrho(p,q)}\|f\|_{M_{p,q}^{[j]}}.   \label{3invariantmod1}
\end{align}
By \eqref{realsubspace} and \eqref{2realsubspace}, to guarantee \eqref{3invariantmod1} is nontrivial, we need to assume that the following condition holds:
\begin{align} \label{3realsubspace}
d\left(\frac{1}{p} + \frac{1}{q} -1\right)    \leq \varrho(p,q) \leq 0.
\end{align}

\begin{figure}
\centering
\includegraphics[width=2.88in, keepaspectratio]{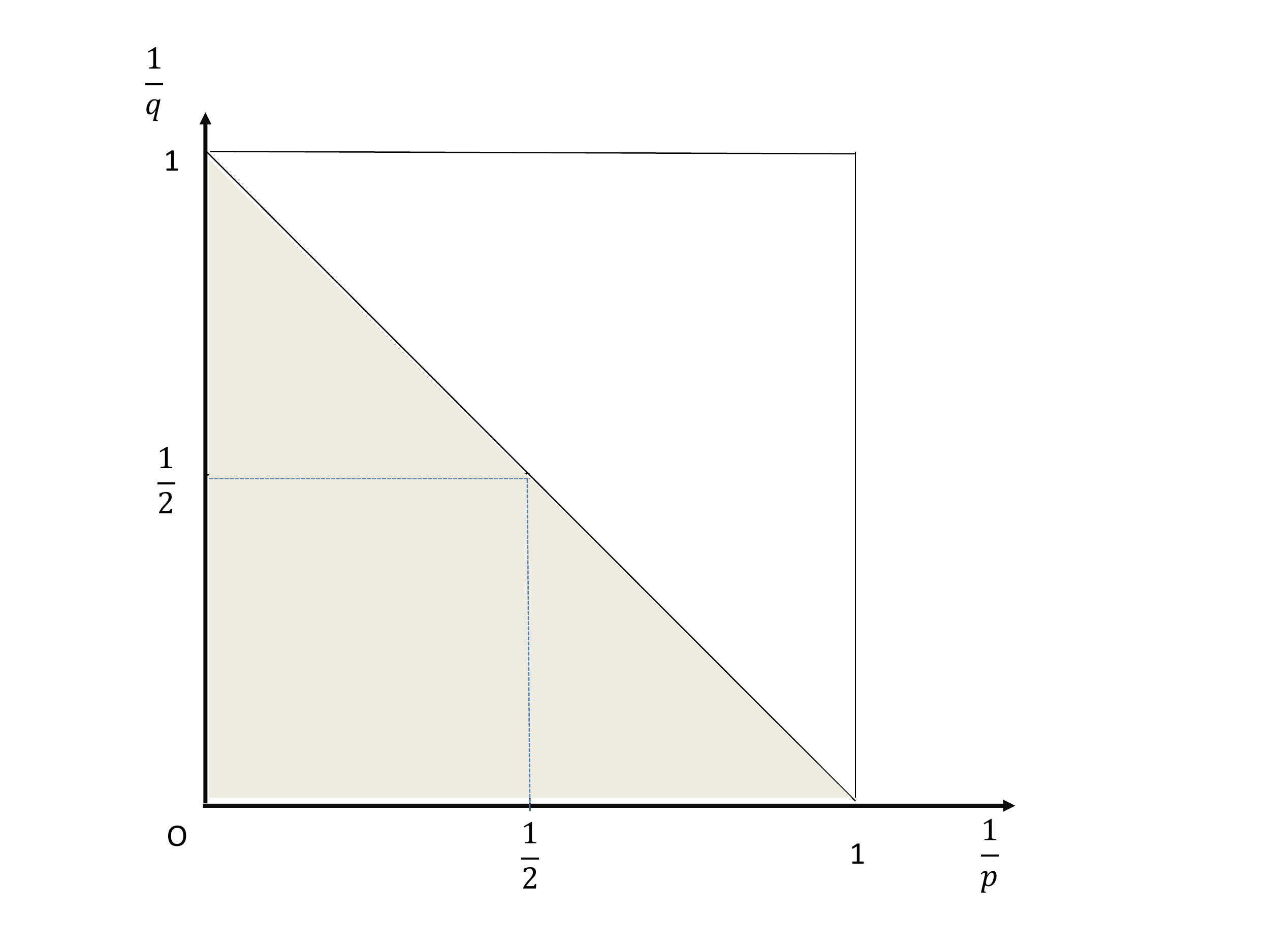}\\
\caption{The grey region of $(1/p,1/q)$: $\varrho(p,q)$ satisfies \eqref{3realsubspace}}
\end{figure}

 In the case $d(1/p + 1/q -1)    < \varrho(p,q) < 0$,  we can introduce a nontrivial subspace $\dot{\mathfrak{M}}_{p,q, 0}^{\varrho}$ of $\dot{\mathfrak{M}}_{p,q, \infty}^{\varrho}$  for which the functions $f \in \dot{\mathfrak{M}}_{p,q, \infty}^{\varrho}$ satisfy the following condition:
\begin{align}
\lim_{j\to \pm \infty} 2^{j\varrho(p,q)}\|f\|_{M_{p,q}^{[j]}}=0.    \label{3invariantmod100}
\end{align}

Combining $ \mathscr{M}^{\mu}_{p,q} $ and $\tilde{\mathscr{M}}^{\mu}_{p,q}$, we assume that $\mu:=\mu(p,q)$ satisfies
\begin{align} \label{3dualcondition}
0 \leq \mu(p,q) \leq   d\left(\frac{1}{q} +\frac{1}{p}-1\right)
\end{align}

\begin{figure}
\centering
\includegraphics[width=2.88in, keepaspectratio]{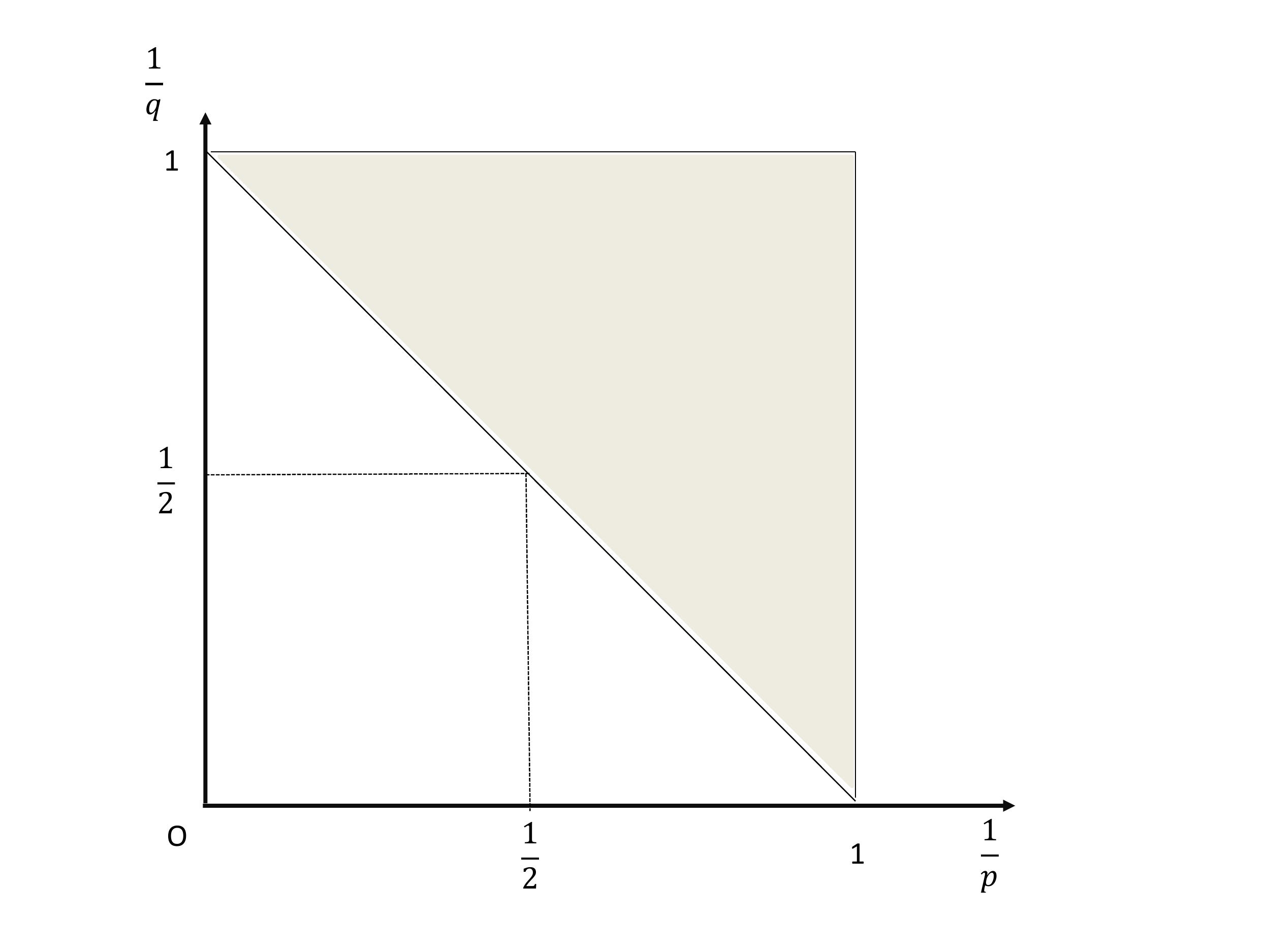}\\
\caption{The grey region of $(1/p,1/q)$ for which $\varrho(p,q)$ satisfies \eqref{3dualcondition}}
\end{figure}

and introduce the following
\begin{align}
 \dot{\mathscr{M}}^{\mu}_{p,q}  = \left\{ f\in \mathscr{S}': \exists \ f_j \in M^{[j]}_{p,q} \ \mbox{such that} \ f= \sum_{j \in \mathbb{Z}} f_j, \ \  \sum_{j \in \mathbb{Z}}    2^{j  \mu } \|f_j\|_{M^{[j]}_{p,q}}<\infty \right\}  \label{3equivinf3}
\end{align}
and the norm on $ \dot{\mathscr{M}}^{\mu}_{p,q}$ is defined as
\begin{align}
\|f\|_{\dot{\mathscr{M}}^{\mu}_{p,q}} = \inf   \sum_{j \in \mathbb{Z}}    2^{j  \mu } \|f_j\|_{M^{[j]}_{p,q}},   \label{3equivinf4}
\end{align}
where the infimum is taken over all of the possible decompositions of $f= \sum_{j\in \mathbb{Z}} f_j \in \dot{\mathscr{M}}^{\mu}_{p,q}$.

\begin{prop} \label{3Banachdense2}
Let $1\leq p,q \leq \infty$, $\varrho(p,q)$ satisfies condition \eqref{3realsubspace}.
Then $\dot{\mathfrak{M}}^{\varrho}_{p,q, \infty}$ and  $\dot{\mathfrak{M}}^{\varrho}_{p,q, 0}$ are Banach spaces.
\end{prop}

\begin{prop} \label{4Banachdense2}
Let $1\leq p,q \leq \infty$, $\mu(p,q)$ satisfies condition \eqref{3dualcondition}.
Then $\dot{\mathscr{M}}^{\mu}_{p,q}$ is a Banach space. Moreover, $\mathscr{S}$ is dense in $\dot{\mathscr{M}}^{\mu}_{p,q}$ if $p,q<\infty$.
\end{prop}

\begin{rem} \rm There is no essential difficulty to generalize the weight  $2^{j\varrho}$ in the definition of ${\mathfrak{M}}^{\varrho}_{p,q,\infty}$ to a general weight ${\bf w} =\{w_j\}$ like
$$
0< w_j \leq
C 2^{jd(1/p+1/q-1)},  \ \  j\leq 0.
$$
By substituting $2^{j\varrho}$ by $w_j$ in the definition of ${\mathfrak{M}}^{\varrho}_{p,q,\infty}$, one get a more general space ${\mathfrak{M}}^{\bf w}_{p,q,\infty}$.  Similarly for $\tilde{{\mathfrak{M}}}^{\bf w}_{p,q,\infty}$ by assuming that
 $$
0< w_j \leq C,  \ \  j\geq 0.
$$
\end{rem}

\section{The duality between $\mathscr{M}^{\mu}_{p,q}$ and $\mathfrak{M}^{-\mu}_{p',q',\infty}$} \label{Duality}
Our main goal of this section is to show the duality between $\mathscr{M}^{\mu}_{p,q}$ and $\mathfrak{M}^{-\mu}_{p',q',\infty}$.
It is easy to see that \eqref{dualcondition} is equivalent to
\begin{align} \label{dualcondition1}
0\vee d\left(\frac{1}{q'} -\frac{1}{p'}\right) \geq  -\mu(p,q)  \geq d\left(\frac{1}{q'} +\frac{1}{p'}-1\right).
\end{align}

\begin{prop} Let $1\leq p,q<\infty$. Assume that $\mu (p,q)$ satisfies \eqref{dualcondition}.  Then we have $(\mathscr{M}^{\mu}_{p,q})^* = \mathfrak{M}^{-\mu}_{p',q',\infty}$ with equivalent norm.
\end{prop}
{\it Proof. } First, we show that  $(\mathscr{M}^{\mu}_{p,q})^* \supset \mathfrak{M}^{-\mu}_{p',q',\infty}$. Let $g\in \mathfrak{M}^{-\mu}_{p',q',\infty}$ and $f= \sum_{j\leq 0} f_j \in \mathscr{M}^{\mu}_{p,q}$. By $\sum_{k\in \mathbb{Z}^d} \Box_{j,k} = I$, we have
\begin{align} \label{dualarg0}
\langle g, f \rangle =\sum_{j\leq 0} \langle g, f_j \rangle =  \sum_{j\leq 0} \sum_{k\in \mathbb{Z}^d} \langle g, \Box_{j,k}f_j \rangle.
\end{align}
Using the almost orthogonal property of $\Box_{j,k}$, and then applying H\"older's inequality, we see that
\begin{align}
\left|\sum_{k\in \mathbb{Z}^d} \langle g, \Box_{j,k}f_j \rangle \right| & \leq  \sum_{k\in \mathbb{Z}^d} \sum_{|\ell|_\infty \leq 1} |\langle \Box_{j,k+\ell}  g, \Box_{j,k}f_j \rangle| \nonumber\\
& \leq  \sum_{k\in \mathbb{Z}^d} \sum_{|\ell|_\infty \leq 1} \|\Box_{j,k+\ell}   g\|_{p'} \|\Box_{j,k}f_j \|_p \nonumber\\
& \lesssim   \| g\|_{M^{[j]}_{p',q'}} \| f_j \|_{M^{[j]}_{p,q}}.  \label{dualarg1}
\end{align}
Hence,
\begin{align}  \label{dualarg2}
|\langle g, f \rangle | \lesssim \sum_{j\leq 0}
   \| g\|_{M^{[j]}_{p',q'}} \| f_j \|_{M^{[j]}_{p,q}} \leq
   \| g\|_{\mathfrak{M}^{-\mu}_{p',q',\infty} }  \sum_{j\leq 0} 2^{j\mu(p,q)}\| f_j \|_{M^{[j]}_{p,q}}.
\end{align}
Taking the infimum over all of the decompositions of $f$, we see that
\begin{align}
|\langle g, f \rangle | \lesssim
   \| g\|_{\mathfrak{M}^{-\mu}_{p',q',\infty} }   \| f\|_{\mathscr{M}^{\mu}_{p,q}}.
\end{align}
It follows that $g\in (\mathscr{M}^{\mu}_{p,q})^*$  and $ \|g\|_{(\mathscr{M}^{\mu}_{p,q})^*} \lesssim \| g\|_{\mathfrak{M}^{-\mu}_{p',q',\infty} }$.

Conversely, for any $g\in (\mathscr{M}^{\mu}_{p,q})^*$, we have
\begin{align}
|\langle g, f \rangle | & \leq
   \| g\|_{ (\mathscr{M}^{\mu}_{p,q})^* }   \| f\|_{\mathscr{M}^{\mu}_{p,q}} \nonumber\\
    & \leq \| g\|_{ (\mathscr{M}^{\mu}_{p,q})^* }  \sum_{j\leq 0}  2^{j\mu(p,q)} \| f_j\|_{M^{[j]}_{p,q}}
\end{align}
for all $f=\sum_{j\leq 0} f_j \in \mathscr{M}^{\mu}_{p,q}$.  Taking $f_j =f$ and $f_\ell =0$ for any $\ell\neq j$, we immediately have
\begin{align}
|\langle g, f \rangle |
    & \leq 2^{j\mu(p,q)} \| g\|_{ (\mathscr{M}^{\mu}_{p,q})^* }   \| f \|_{M^{[j]}_{p,q}}
\end{align}
Using the duality $(M^{[j]}_{p,q})^* = M^{[j]}_{p',q'}$ (see Appendix), we have
\begin{align}
 2^{-j\mu(p,q)}\| g \|_{M^{[j]}_{p',q'}} \lesssim   \| g\|_{ (\mathscr{M}^{\mu}_{p,q})^* }.
\end{align}
Hence, we have $g\in \mathfrak{M}^{-\mu}_{p',q',\infty}$ and $\|g\|_{\mathfrak{M}^{-\mu}_{p',q',\infty}} \lesssim   \| g\|_{ (\mathscr{M}^{\mu}_{p,q})^* }$.

\begin{rem} \rm
Recall that the pre-dual of Fofana spaces were first considered in \cite{FeFo19} by using the ideas of exotic minimal space developed in \cite{FeZi02}. Kpata \cite{Kp20} considered the predual of $(L^1, \ell^q)^\alpha$ by following  \cite{FeFo19}. The techniques in \cite{FeZi02} seem also useful for the scaling limit spaces of $M^0_{p,q}$. On the other hand, the results in this Section can be developed to the  Fofana spaces without essential difficulty.

\end{rem}

\begin{prop}
Let $1\leq p,q <\infty$ and $\varrho(p,q)$ satisfy \eqref{realsubspace} and $\varrho\neq d(1/p+1/q-1)$. Then we have
$(\mathfrak{M}^{\varrho}_{p,q,0})^* =  \mathscr{M}^{-\varrho}_{p',q'}$.
\end{prop}
{\it Proof.} By the dual argument as in \eqref{dualarg1} and  \eqref{dualarg2}, we see that $\mathscr{M}^{-\varrho}_{p',q'} \subset (\mathfrak{M}^{\varrho}_{p,q,0})^*$ and
$$
\|g\|_{(\mathfrak{M}^{\varrho}_{p,q,0})^*} \lesssim  \|g\|_{\mathscr{M}^{-\varrho}_{p',q'}}, \ \ \forall g\in \mathscr{M}^{-\varrho}_{p',q'}.
$$
Conversely, let us write $\ell^{\infty, \varrho}_0(\mathbb{Z}_-; \ell^q(\mathbb{Z}^d, L^p(\mathbb{R}^d)))$ as the space of all function sequences $f:= \{f_{j,k}\}$ satisfying
$$
\|f\|_{\ell^{\infty, \varrho}_0(\mathbb{Z}_-; \ell^q(\mathbb{Z}^d, L^p(\mathbb{R}^d)))} = \sup_{j\leq 0} 2^{j\varrho}
\left(\sum_{k\in\mathbb{Z}^d} \|f_{j,k}\|^q_{L^p(\mathbb{R}^d)}\right)^{1/q} <\infty
$$
and
$$
\lim_{j\to -\infty} \left(\sum_{k\in \mathbb{Z}^d} \|f_{j,k}\|^q_{L^p(\mathbb{R}^d)}\right)^{1/q}=0.
$$
One can treat $\mathfrak{M}^{\varrho}_{p,q,0}$ as a subspace of $\ell^{\infty, \varrho}_0(\mathbb{Z}_-; \ell^q(\mathbb{Z}^d, L^p(\mathbb{R}^d)))$ via the isometric mapping:
$$
\Gamma: \mathfrak{M}^{\varrho}_{p,q,0} \ni f \mapsto \{\Box_{j,k} f\}_{j\leq 0, k\in \mathbb{Z}^d} \in  \ell^{\infty, \varrho}_0(\mathbb{Z}_-; \ell^q(\mathbb{Z}^d, L^p(\mathbb{R}^d))).
$$
So, $g\in (\mathfrak{M}^{\varrho}_{p,q,0})^*$ can be regarded as a continuous functional on
$$
X_0 = \{\{\Box_{j,k} f\}_{j\leq 0, k\in \mathbb{Z}^d}: f\in \mathfrak{M}^{\varrho}_{p,q,0}\}.
$$
By Hahn-Banach's Theorem, $g$  can be extended onto $\ell^{\infty, \varrho}_0(\mathbb{Z}_-; \ell^q(\mathbb{Z}^d, L^p(\mathbb{R}^d)))$ and the extension is written as $\tilde{g}$, i.e., $\tilde{g} \in (\ell^{\infty, \varrho}_0(\mathbb{Z}_-; \ell^q(\mathbb{Z}^d, L^p(\mathbb{R}^d))))^* = \ell^{1, -\varrho}_0(\mathbb{Z}_-; \ell^{q'}(\mathbb{Z}^d, L^{p'}(\mathbb{R}^d)))$, where $\tilde{g} = \{g_{j,k}\}_{j\leq 0, k\in \mathbb{Z}^d}$ is equipped with the norm
$$
\|\tilde{g}\|_{\ell^{1, -\varrho}_0(\mathbb{Z}_-; \ell^{q'}(\mathbb{Z}^d, L^{p'}(\mathbb{R}^d)))} = \sum_{j\leq 0} 2^{-j\varrho}
\left(\sum_{k\in \mathbb{Z}^d} \|g_{j,k}\|^{q'}_{L^{p'}(\mathbb{R}^d)}\right)^{1/q'},
$$
and $\|g\|_{(\mathfrak{M}^{\varrho}_{p,q,0})^*} = \|\tilde{g}\|_{\ell^{1, -\varrho}_0(\mathbb{Z}_-; \ell^{q'}(\mathbb{Z}^d, L^{p'}(\mathbb{R}^d)))}$. Moreover,
$$
\langle \tilde{g}, \{f_{j,k}\}\rangle = \sum_{j\leq 0}\sum_{k\in \mathbb{Z}^d} \int_{\mathbb{R}^d} f_{j,k}(x) \overline{g_{j,k}(x)} dx.
$$
So, if we restrict $\tilde{g}$ on $X_0$, we have for any $f\in \mathfrak{M}^{\varrho}_{p,q,0}$,
\begin{align}
\langle g, f\rangle & = \langle \tilde{g}, \{\Box_{j,k}f\}\rangle = \sum_{j\leq 0}\sum_{k\in \mathbb{Z}^d} \int_{\mathbb{R}^d} \Box_{j,k}f (x) \overline{g_{j,k}(x)} dx.  \nonumber
\end{align}
It follows that for any $f\in \mathscr{S}$,
\begin{align}
\langle g, f\rangle =  \sum_{j\leq 0} \langle \sum_{k\in \mathbb{Z}^d}\overline{\Box_{j,k} g_{j,k} }, \  f \rangle.
\end{align}
Let us write $g_j =  \sum_{k\in \mathbb{Z}^d} \Box_{j,k} g_{j,k}$, we see that $g= \sum_{j\leq 0} g_j$ and
$$
\|g_j\|_{M^{[j]}_{p',q'}} \leq  \left(\sum_{k\in \mathbb{Z}^d} \|g_{j,k}\|^{q'}_{L^{p'}(\mathbb{R}^d)}\right)^{1/q'},
$$
which implies that
$$
\sum_{j\leq 0} 2^{-j\varrho}\|g_j\|_{M^{[j]}_{p',q'}} \leq  \|\tilde{g}\|_{\ell^{1, -\varrho}_0(\mathbb{Z}_-; \ell^{q'}(\mathbb{Z}^d, L^{p'}(\mathbb{R}^d)))} = \|g\|_{(\mathfrak{M}^{\varrho}_{p,q,0})^*}.
$$
Hence, $g\in \mathscr{M}^{-\varrho}_{p',q'}$ and $\|g\|_{\mathscr{M}^{-\varrho}_{p',q'}} \lesssim \|g\|_{(\mathfrak{M}^{\varrho}_{p,q,0})^*}.$ $\hfill\Box$ \\

\begin{cor}
Let $1\leq p,q <\infty$ and $\varrho(p,q)$ satisfy \eqref{realsubspace}. Then we have
$(\mathfrak{M}^{\varrho}_{p,q,\infty})^* =  \mathscr{M}^{-\varrho}_{p',q'} + (\mathfrak{M}^{\varrho}_{p,q,0})^\bot$, where  $$(\mathfrak{M}^{\varrho}_{p,q,0})^\bot = \{f\in  (\mathfrak{M}^{\varrho}_{p,q,\infty})^*: \, \langle f, \, \varphi\rangle =0 , \ \forall \varphi \in \mathfrak{M}^{\varrho}_{p,q,0}\}.$$
\end{cor}
{\bf Proof.} Let $X$ be a Banach space, $L$ is a subspace of $X$,
$$L^\bot = \{f\in X^*:\, \langle f, \, \phi\rangle =0, \ \forall \phi \in L\}.$$
By a well-known dual result, we see that $L^* = X^*/ L^\bot$. This implies our result. $\hfill\Box$\\

Following the arguments as above, we can show the following results and the proofs will be omitted.

\begin{prop} Let $1\leq p,q<\infty$. Assume that $\mu (p,q)$ satisfies \eqref{2dualcondition}.  Then we have $(\tilde{\mathscr{M}}^{\mu}_{p,q})^* = \tilde{\mathfrak{M}}^{-\mu}_{p',q',\infty}$ with equivalent norms.
\end{prop}

\begin{prop} Let $1\leq p,q<\infty$. Assume that $\mu (p,q)$ satisfies \eqref{3dualcondition}.  Then we have $(\dot{\mathscr{M}}^{\mu}_{p,q})^* = \dot{\mathfrak{M}}^{-\mu}_{p',q',\infty}$ with equivalent norms.
\end{prop}

\begin{prop}
Let $1\leq p,q <\infty$ and $\varrho(p,q)$ satisfy \eqref{2realsubspace} and $\varrho\neq 0$. Then we have
$(\tilde{\mathfrak{M}}^{\varrho}_{p,q,0})^* =  \tilde{\mathscr{M}}^{-\varrho}_{p',q'}$ with equivalent norms.
\end{prop}

\begin{prop}
Let $1\leq p,q <\infty$ and $\varrho(p,q)$ satisfy $d(1/p+1/q-1) < \varrho (p,q) <0$. Then we have
$(\dot{\mathfrak{M}}^{\varrho}_{p,q,0})^* =  \dot{\mathscr{M}}^{-\varrho}_{p',q'}$ with equivalent norms.
\end{prop}
In a similar way as Corollary 3.4, we can describe the dual spaces of $\tilde{\mathfrak{M}}^{\varrho}_{p,q,\infty}$ and $\dot{\mathfrak{M}}^{\varrho}_{p,q,\infty}$ and the details are omitted.

\section{Generalizations of $\mathscr{M}^{\mu}_{p,q}$ and $\mathfrak{M}^{\varrho}_{p,q, \infty}$}
\label{Generalmodspace}

We consider the generalization of  $\mathscr{M}^{\mu}_{p,q}$ and introduce the following
\begin{align}
 \mathscr{M}^{\mu}_{p,q, r}  = \left\{ f\in \mathscr{S}': \exists \ f_j \in M^{[j]}_{p,q} \ \mbox{such that} \ f= \sum_{j\leq 0} f_j, \ \  \sum_{j\leq 0}    2^{j r \mu } \|f_j\|^r_{M^{[j]}_{p,q}}<\infty \right\}  \label{4equivinf3}
\end{align}
and the norm on $ \mathscr{M}^{\mu}_{p,q, r}$ is defined as
\begin{align}
\|f\|_{\mathscr{M}^{\mu}_{p,q, r}} = \inf   \left(\sum_{j\leq 0}    2^{j r \mu } \|f_j\|^r_{M^{[j]}_{p,q}} \right)^{1/r},   \label{4equivinf4}
\end{align}
where the infimum is taken over all of the decompositions of $f= \sum_{j\leq 0} f_j \in \mathscr{M}^{\mu}_{p,q,r}$.
It is easy to see that $ \mathscr{M}^{\mu}_{p,q}=  \mathscr{M}^{\mu}_{p,q, 1}$.  Following B\'{e}nyi and Oh's \cite{BeOh20}, we define
\begin{align}
\|f\|_{\mathfrak{M}^{\varrho}_{p,q, r}} =    \left(\sum_{j\leq 0}    2^{j r \varrho } \|f\|^r_{M^{[j]}_{p,q}} \right)^{1/r}.     \label{4equivinf5}
\end{align}

\begin{prop} \label{Embedding1}
 Let $1\leq  r<\infty$, $1\leq p,q\leq \infty$. Assume that $\mu (p,q)$ satisfies
 \begin{align} \label{4embedding0}
 0\wedge d\left(\frac{1}{q}- \frac{1}{p}\right) < \mu(p,q) <  d\left(\frac{1}{q}+ \frac{1}{p} -1\right).
 \end{align}
 Then we have for any $\tilde{p} \geq p$  with $\mu(p,q) <  d (1/p+1/q -1/\tilde{p} -1 ) $,
 \begin{align} \label{4embedding1}
 M^{0}_{p,q} \subset  \mathscr{M}^{\mu}_{p,q, r}  \subset M^0_{\tilde{p}, q}.
 \end{align}
 Moreover, \eqref{4embedding1} also hold for  $\mu(p,q) = d (1/p+1/q -1/\tilde{p} -1)$ if $r=1$.
\end{prop}

Proposition \ref{Embedding1} indicates the relations between  $ \mathscr{M}^{\mu}_{p,q, r}$ and  $ M^{0}_{p,q}$. In particular, \eqref{4embedding1} implies that the following embeddings hold:
 \begin{align} \label{4embedding2}
 M^{0}_{p,q} \subset  \mathscr{M}^{\mu}_{p,q, r}  \subset M^0_{\infty, q}.
 \end{align}

{\it Proof of Proposition \ref{Embedding1}.} Denote
$$
\Lambda_{j,k} = \{\ell \in \mathbb{Z}^d: \ 2^j(\ell + [-1/2, 1/2]^d) \cap (k + [-1, 1]^d)\neq \emptyset\}.
$$
It follows from Bernstein's inequality that
$$
\|\Box_k g\|_{\tilde{p}} \leq \sum_{\ell \in \Lambda_{j,k}} \|\Box_{j,\ell}g\|_{\tilde{p}}  \lesssim 2^{jd(1/p-1/\tilde{p})} \sum_{\ell \in \Lambda_{j,k}} \|\Box_{j,k}f\|_{p}.
$$
By H\"older's inequality, we have
$$
\|\Box_k g\|_{\tilde{p}}   \lesssim 2^{jd(1/p +1/q -1/\tilde{p}-1)} \left( \sum_{\ell \in \Lambda_{j,k}} \|\Box_{j,k}g\|^q_{p} \right)^{1/q}.
$$
Taking the $\ell^q$ norm in both sides of the above inequality, we have
$$
\|g\|_{M^{0}_{\tilde{p},q}} \lesssim 2^{jd(1/p +1/q -1/\tilde{p}-1)} \left( \sum_{k\in \mathbb{Z}^d} \sum_{\ell \in \Lambda_{j,k}} \|\Box_{j,k}g\|^q_{p} \right)^{1/q}
$$
Since $\Lambda_{j,k} \cap \Lambda_{j,k+l} \neq \emptyset$ implies that $|\ell|_\infty \lesssim 1$, we have
\begin{align} \label{4embedding3}
\|g\|_{M^{0}_{\tilde{p},q}} \lesssim 2^{jd(1/p +1/q -1/\tilde{p}-1)} \|g \|_{M^{[j]}_{p,q}}.
\end{align}
We emphasize that the above estimates uniformly hold for all $j\leq 0$. For any $f = \sum_{j\leq 0} f_j \in \mathscr{M}^{\mu}_{p,q,r}$, using H\"older's inequality, we have from \eqref{4embedding3} that
\begin{align}
\|f\|_{M^{0}_{\tilde{p},q}}  & \leq \sum_{j\leq 0 } \|f_j\|_{M^{0}_{\tilde{p},q}}  \nonumber\\
& \lesssim \sum_{j\leq 0 } 2^{jd(1/p +1/q -1/\tilde{p}-1)} \|f_j\|_{M^{[j]}_{p,q}} \nonumber\\
& \lesssim    \left\|2^{j\mu(p,q)}\|f_j\|_{M^{[j]}_{p,q}} \right\|_{\ell^r_j},  \label{4embedding4}
\end{align}
where we used the fact that $d(1/p +1/q -1/\tilde{p}-1) -\mu(p,q)>0$ and \eqref{4embedding4} also hold for $r=1$ if $d(1/p +1/q -1/\tilde{p}-1) =\mu(p,q)$.  Taking the infimum on all of the possible decomposition of $f \in \mathscr{M}^{\mu}_{p,q,r}$ in \eqref{4embedding4}, we have
$$
\|f\|_{M^{0}_{\tilde{p},q}}  \lesssim    \|f\|_{\mathscr{M}^{\mu}_{p,q,r}}.
$$
Hence, $\mathscr{M}^{\mu}_{p,q,r} \subset M^{0}_{\tilde{p},q}$. The embedding $M^{0}_{p,q} \subset  \mathscr{M}^{\mu}_{p,q,r}$ is obtained by considering the decomposition of $f=f +0+...+0+...$ .  $\hfill\Box$

\begin{prop}
Let $1\leq p,q \leq \infty$, $1<r<\infty$. Assume that $\varrho: = \varrho(p,q)$ and $\mu:=\mu (p,q)$  satisfy \eqref{realsubspace} and \eqref{dualcondition}, respectively. $\varrho \neq d(1/p+1/q-1)$. Then  $\mathfrak{M}^{\varrho}_{p,q,r}$ and $ \mathscr{M}^{\mu}_{p,q, r} $ are Banach spaces. Moreover, $\mathscr{S}$ is dense in $\mathfrak{M}^{\varrho}_{p,q,r}$ and $ \mathscr{M}^{\mu}_{p,q, r} $ if $p,q<\infty$.
\end{prop}
{\bf Proof.} Following the same way as in the proofs of  Propositions \ref{Banach-1} and \ref{Banachdense2}, we have the results, as desired.  $\hfill\Box$

\begin{prop}
Let $1< p,q, r<\infty$. Assume that $\mu:=\mu (p,q)$ satisfies \eqref{dualcondition} and $\mu < d(1/p+1/q-1)$.  Then we have $(\mathscr{M}^{\mu}_{p,q, r})^* = \mathfrak{M}^{-\mu}_{p',q',r'}$ and  $(\mathfrak{M}^{-\mu}_{p',q',r'})^* = \mathscr{M}^{\mu}_{p,q, r}  $ with equivalent norm.
\end{prop}
{\it Sketch Proof. } For any $g \in \mathfrak{M}^{-\mu}_{p',q',r'}$ and $f=\sum_{j\leq 0} f_j \in \mathscr{M}^{\mu}_{p,q, r}$,  in view of \eqref{dualarg0},\eqref{dualarg1} and  H\"older's inequality,
\begin{align}  \label{4dualarg2}
|\langle g, f \rangle | \lesssim \sum_{j\leq 0}
   \| g\|_{M^{[j]}_{p',q'}} \| f_j \|_{M^{[j]}_{p,q}} \leq \|g\|_{\mathfrak{M}^{-\mu}_{p',q',r'}}  \left\| 2^{j\mu(p,q)} \| f_j \|_{M^{[j]}_{p,q}} \right\|_{\ell^r_j}.
\end{align}
Taking the infimum on all of the possible decomposition of $f \in \mathscr{M}^{\mu}_{p,q,r}$, by \eqref{4dualarg2} we have $\mathfrak{M}^{-\mu}_{p',q',r'} \subset (\mathscr{M}^{\mu}_{p,q, r})^*$.  On the other hand, for any $g\in (\mathscr{M}^{\mu}_{p,q, r})^*$, we have
\begin{align} \label{4dualarg1}
|\langle g, f \rangle | \leq \|g\|_{ (\mathscr{M}^{\mu}_{p,q, r})^* } \left\| 2^{j\mu(p,q)} \| f_j \|_{M^{[j]}_{p,q}} \right\|_{\ell^r_j}
\end{align}
for all $f = \sum_{j\leq 0} f_j  \in  \mathscr{M}^{\mu}_{p,q, r}$.
By Proposition \ref{Embedding1}, we see that $(\mathscr{M}^{\mu}_{p,q, r})^* \subset M^0_{p'q'}$. We have for any $\varphi=\sum_{-J\leq j\leq 0} \varphi_j, \ \varphi_j \in \mathscr{S}$,
\begin{align} \label{4dualarg3}
\left|\sum^0_{j=-J} \langle g, \varphi_j  \rangle  \right|  \leq \|g\|_{ (\mathscr{M}^{\mu}_{p,q, r})^* } \left(\sum^0_{j=-J} 2^{j r\mu(p,q)} \| \varphi_j \|^r_{M^{[j]}_{p,q}} \right)^{1/r}.
\end{align}
Since $\mathscr{S}$ is dense in  $ {M}^{[j]}_{p,q}$, by duality we have,
\begin{align} \label{4dualarg4}
\left(\sum^0_{j=-J} 2^{- j r' \mu(p,q)} \left\|  g \right\|^{r'}_{M^{[j]}_{p',q'}} \right)^{1/r'}     \leq  \|g\|_{ (\mathscr{M}^{\mu}_{p,q, r})^* }.
\end{align}
Letting $J\to \infty$ in \eqref{4dualarg4}, we immediately have
$$
 \|g\|_{ \mathscr{M}^{-\mu}_{p',q', r'} }  \leq  \|g\|_{ (\mathscr{M}^{\mu}_{p,q, r})^* }.
$$
The result follows. $\hfill\Box$

Following the arguments as above, we can show the following results and the proofs will be omitted.

\begin{prop} Let $1< p,q,r<\infty$. Assume that $\mu:=\mu (p,q)$ satisfies \eqref{2dualcondition}.  Then we have $(\tilde{\mathscr{M}}^{\mu}_{p,q,r})^* = \tilde{\mathfrak{M}}^{-\mu}_{p',q', r'}$ with equivalent norms.
\end{prop}

\begin{prop} Let $1< p,q, r <\infty$. Assume that $\mu:=\mu (p,q)$ satisfies \eqref{3dualcondition}.  Then we have $(\dot{\mathscr{M}}^{\mu}_{p,q,r})^* = \dot{\mathfrak{M}}^{-\mu}_{p',q', r}$ with equivalent norms.
\end{prop}

\begin{prop}
Let $1< p,q, r <\infty$ and $\varrho:=\varrho(p,q)$ satisfy \eqref{2realsubspace} and $\varrho \neq 0$. Then we have
$(\tilde{\mathfrak{M}}^{\varrho}_{p,q,r})^* =  \tilde{\mathscr{M}}^{-\varrho}_{p',q',r'}$ with equivalent norms.
\end{prop}

\begin{prop}
Let $1< p,q <\infty$ and $\varrho:=\varrho(p,q)$ satisfy $d(1/p+1/q-1) < \varrho <0$. Then we have
$(\dot{\mathfrak{M}}^{\varrho}_{p,q,r})^* =  \dot{\mathscr{M}}^{-\varrho}_{p',q',r'}$ with equivalent norms.
\end{prop}

\begin{rem} \rm There is no essential difficulty to generalize the weight  $2^{j\varrho}$ in the definition of ${\mathfrak{M}}^{\varrho}_{p,q,r}$ to a general weight ${\bf w} =\{w_j\}$ like
$$
0< w_j \leq
C 2^{jd(1/p+1/q-1 + \varepsilon)},  \ \  j\leq 0.
$$
By substituting $2^{j\varrho}$ by $w_j$ in the definition of ${\mathfrak{M}}^{\varrho}_{p,q,r}$, one get a more general space ${\mathfrak{M}}^{\bf w}_{p,q,r}$.  Similarly for $\tilde{{\mathfrak{M}}}^{\bf w}_{p,q,r}$ by assuming that
 $$
0< w_j \leq C 2^{- \varepsilon j},  \ \  j\geq 0.
$$
One can get that their dual spaces  are $ {\mathscr{M}}^{ \{1/w_j\}}_{p',q',r'}$  and  $\tilde{{\mathscr{M}}}^{ \{1/w_j\}}_{p',q',r'}$ for a class general weights and we will not perform this generalization in this paper.
\end{rem}

\section{Scalings} \label{Scaling}

Let us start with a useful equivalent norm in  $\mathscr{M}^{\mu}_{p,q, r}$ and $\mathfrak{M}^{\rho}_{p,q, r}$.  Denote
\begin{align}
\|f\|^{c}_{M^{[j]}_{p,q}} = \left\| \|\mathscr{F}^{-1}\psi(2^{-j} c \ \cdot -k ) \hat{f}\|_p \right\|_{\ell^q}.
\end{align}

For any $c>0$, it is easy to see that ${\rm supp} \psi(2^{-j} c \ \cdot -k )$ overlaps at most finite many ${\rm supp} \psi_{j,l}$  and vice versa. From Bernstein's multiplier estimate, we see that $\|\cdot\|^{c}_{M^{[j]}_{p,q}}  $ and $ \|\cdot\|_{M^{[j]}_{p,q}}$ are equivalent norms and
\begin{align} \label{equivdia}
c_0 \|f\|_{M^{[j]}_{p,q}} \leq \|f\|^{c}_{M^{[j]}_{p,q}} \leq C_1 \|f\|_{M^{[j]}_{p,q}}
\end{align}
holds for all $f\in M^{[j]}_{p,q}$, $c_0$ and $C_1$ are independent of $j\in \mathbb{Z}$ and $c\in [1,2]$. So, we immediately have

\begin{lem}
Let $c>0$.   Then
\begin{align}
\|f\|^c_{\mathfrak{M}^{\varrho}_{p,q, r}} = \left \| 2^{j\varrho}  \|f\|^{c}_{M^{[j]}_{p,q}}  \right\|_{\ell^r_j}
\end{align}
is an equivalent norm on $\mathfrak{M}^{\varrho}_{p,q, r}$. The same result holds on  $\mathscr{M}^{\mu}_{p,q, r}$ if one replaces $\|\cdot\|_{M^{[j]}_{p,q}}$ with $\|\cdot\|^{c}_{M^{[j]}_{p,q}}$ in the definition of $\mathscr{M}^{\mu}_{p,q, r}$.
\end{lem}
We have the following scaling property of  $\mathfrak{M}^{\varrho}_{p,q, r}$.

\begin{prop} \label{ScalingP1}
Let $1\leq p,q\leq \infty$, $1 < r \leq \infty$ and $\varrho:= \varrho(p,q)$ satisfy \eqref{realsubspace}, and $\varrho\neq d(1/p+1/q-1)$ for $1<r<\infty$. Then we have for any $\lambda\geq 1$,
$$
\|f_\lambda\|_{\mathfrak{M}^{\varrho}_{p,q, r}} \leq C \lambda^{\varrho -d/p} \|f \|_{\mathfrak{M}^{\varrho}_{p,q, r}},
$$
where $C>1$ is independent of $\lambda\geq 1$.
\end{prop}
{\it Proof.} Let $\lambda = 2^{-j_0} c$ with $j_0 \in -\mathbb{Z}_-$ and $c \in [1,2)$. Similar to \eqref{dialambda},
\begin{align} \label{scalingdialambda}
\| f_\lambda\|_{M^{[j]}_{p,q}}= \lambda^{-d/p}\left(\sum_{k\in \mathbb{Z}^d}  \|\mathscr{F}^{-1} \psi (\lambda 2^{-j} \xi -k) \mathscr{F} f\|^q_p \right)^{1/q}.
\end{align}
Hence, in view of \eqref{equivdia},
\begin{align} \label{scalingdialambda1}
 \left\| 2^{j\varrho} \| f_\lambda\|_{M^{[j]}_{p,q}} \right\|_{\ell^r_{j\leq 0}} & = \lambda^{-d/p} 2^{-j_0\varrho} \left\| 2^{j \varrho} \|f\|^{c}_{M^{[j]}_{p,q}}  \right\|_{\ell^r_{j\leq j_0}}\nonumber\\
& \leq C \lambda^{\varrho -d/p}   \|f\|_{\mathfrak{M}^{\varrho}_{p,q, r}}.
\end{align}
The result follows. $\hfill\Box$\\

The next two results can be obtained by following the same way as that of Proposition \ref{ScalingP1} and we omit the details of their proofs.

\begin{prop} \label{ScalingP2}
Let $1\leq p,q\leq \infty$, $1 < r \leq \infty$ and $\varrho:= \varrho(p,q)$ satisfy \eqref{2realsubspace}, and $\varrho\neq 0$ for $1<r<\infty$. Then we have for any $0< \lambda\leq 1$,
$$
\|f_\lambda\|_{\tilde{\mathfrak{M}}^{\varrho}_{p,q, r}} \leq C \lambda^{\varrho -d/p} \|f \|_{\tilde{\mathfrak{M}}^{\varrho}_{p,q, r}},
$$
where $C>1$ is independent of $\lambda\leq 1$.
\end{prop}

\begin{prop} \label{ScalingP3}
Let $1\leq p,q\leq \infty$, $1 < r \leq \infty$ and $\varrho$ satisfy $d(1/p+1/q-1)<\varrho <0$. Then we have for any $\lambda>0$,
$$
\|f_\lambda\|_{\dot{\mathfrak{M}}^{\varrho}_{p,q, r}} \sim \lambda^{\varrho -d/p} \|f \|_{\dot{\mathfrak{M}}^{\varrho}_{p,q, r}},
$$
where $C>1$ is independent of $\lambda >0$.
\end{prop}

Using the first inequality in \eqref{scalingdialambda1}, or by the duality, we have

\begin{prop} \label{ScalingP4}
Let $1\leq p,q < \infty$, $1 \leq r < \infty$ and $\mu:=\mu(p,q)$ satisfy \eqref{dualcondition}, \eqref{2dualcondition} \eqref{3dualcondition} for $\mathscr{M}^{\mu}_{p,q, r}$, $\tilde{\mathscr{M}}^{\mu}_{p,q, r}$ and $ \dot{\mathscr{M}}^{\mu}_{p,q, r}$, respectively.  Then we have
\begin{align}
\|f_\lambda\|_{\mathscr{M}^{\mu}_{p,q, r}}  & \lesssim \lambda^{\mu -d/p} \|f \|_{\mathscr{M}^{\mu}_{p,q, r}}, \ \ \forall  \ \lambda >1, \\
\|f_\lambda\|_{\tilde{\mathscr{M}}^{\mu}_{p,q, r}}  & \lesssim \lambda^{\mu -d/p} \|f \|_{\tilde{\mathscr{M}}^{\mu}_{p,q, r}}, \ \ \forall \ 0< \lambda \leq 1. \\
\|f_\lambda\|_{\dot{\mathscr{M}}^{\mu}_{p,q, r}}  & \sim \lambda^{\mu -d/p} \|f \|_{\dot{\mathscr{M}}^{\mu}_{p,q, r}}, \ \ \forall  \ \lambda >0.
\end{align}
\end{prop}

\section{Algebraic property of $\mathscr{M}^{\mu}_{p,q,r}$} \label{AlgebraP}

In this section we consider an algebraic  structure of $\mathscr{M}^{\mu}_{p,q,r}$, which is of importance in the study of nonlinear PDE.

\begin{lem} \label{Algebra1}
Let $1\leq p,q\leq \infty$, $i\leq j\leq 0$. Denote
$$
\mu_0 (p,q) = 0\wedge d \left(\frac{1}{q} -\frac{1}{p}\right) \wedge d \left(\frac{1}{q}+\frac{1}{p}-1\right).
$$
Then we have
$$
\|f\|_{M^{[j]}_{p,q}} \lesssim 2^{-(j-i)\mu_0(p,q)} \|f\|_{M^{[i]}_{p,q}}
$$
for all $g\in M^{[i]}_{p,q}$.
\end{lem}
{\it Proof.} In view of \eqref{inmodulation-j}, we have for $\lambda=2^{-j}$, $\tau =2^{-i}$,
$$
\|f\|_{M^{[j]}_{p,q}} = \lambda^{d/p} \|f_\lambda\|_{M^{0}_{p,q}}, \ \ \ \|f\|_{M^{[i]}_{p,q}} = \tau^{d/p} \|f_\tau\|_{M^{0}_{p,q}}.
$$
Moreover, by the first inequality in the dilation property \eqref{dimo},
$$
\|f_\tau\|_{M^{0}_{p,q}} = \|(f_\lambda)_{2^{j-i}}\|_{M^{0}_{p,q}} \gtrsim 2^{-d(j-i)/p } 2^{(j-i)\mu_0(p,q)} \|f_\lambda\|_{M^{0}_{p,q}}.
$$
Combining the above estimates, we have the result, as desired. $\hfill\Box$

\begin{prop} \label{Algebra2}
Let $1\leq p, r<\infty$. Assume that $0\leq \mu < d/p$ for $r\in (1,\infty)$; $0\leq \mu  \leq d/p$ for $r=1$. Then $\mathscr{M}^{\mu}_{p,1, r}$ is a Banach algebra.
\end{prop}

{\it Proof.}  Let $f,g\in \mathscr{M}^{\mu}_{p,1, r}$ with $f=\sum_{j\leq 0} f_j$ and $g=\sum_{i \leq 0} g_i$.  Further, we can assume that
$$
\left(\sum_{j\leq 0} 2^{j\mu r} \|f_j\|^r_{M^{[j]}_{p,1}} \right)^{1/r} \leq 2 \|f\|_{\mathscr{M}^{\mu}_{p,1,r}}, \ \ \ \left(\sum_{i  \leq 0} 2^{i \mu r} \|g_i\|^r_{M^{[i]}_{p,1}} \right)^{1/r} \leq 2 \|g\|_{\mathscr{M}^{\mu}_{p,1,r}}.
$$
We have
\begin{align}
fg   = \sum_{j\leq 0} \sum_{i\leq j} f_j g_i + \sum_{i\leq 0} \sum_{j<i} f_j g_i
  :=  \sum_{j\leq 0}  F_j +   \sum_{i\leq 0} G_i .  \label{fgdecom}
\end{align}
It follows that
\begin{align}
\|fg\|_{\mathscr{M}^{\mu}_{p,1,r}} &  \leq  \left(\sum_{j\leq 0} 2^{j\mu r} \|F_j\|^r_{M^{[j]}_{p,1}} \right)^{1/r}  +  \left( \sum_{i\leq 0}   2^{i \mu r} \|G_i\|^r_{M^{[i]}_{p,1}} \right)^{1/r}.  \label{fgalgebra1}
\end{align}
One needs to estimate $\|F_j\|_{M^{[j]}_{p,1}}$. We have
\begin{align}
 \|F_j\|_{M^{[j]}_{p,1}}  \leq  \sum_{i\leq j}   \|f_j g_i\|_{M^{[j]}_{p,1}}.  \label{fgalgebra2}
\end{align}
For any $i\leq j\leq 0$, we have
\begin{align}
 \|f_j g_i\|_{M^{[j]}_{p,1}} = \sum_{k\in \mathbb{Z}^d} \|\Box_{j,k} (f_j g_i)\|_p  .  \label{fgalgebra3}
\end{align}
Applying the almost orthogonality of $\Box_{j,k}$, we have
\begin{align}
  \|\Box_{j,k} (f_j g_i)\|_p   \leq  \sum_{|\ell|_\infty \leq 2} \sum_{k_1 \in \mathbb{Z}^d} \|\Box_{j,k} (\Box_{j,k_1}  f_j  \ \Box_{j,k-k_1-\ell}  g_i) \|_p .  \label{fgalgebra4}
\end{align}
Since $\psi_{j,k}$ is a multiplier in $L^p$, we see that
\begin{align}
  \|\Box_{j,k} (f_j g_i)\|_p  &  \leq  \sum_{|\ell|_\infty \leq 2} \sum_{k_1 \in \mathbb{Z}^d} \|  \Box_{j,k_1}  f_j  \ \Box_{j,k-k_1-\ell}  g_i  \|_p \nonumber\\
  & \leq  \sum_{|\ell|_\infty \leq 2} \sum_{k_1 \in \mathbb{Z}^d} \|  \Box_{j,k_1}  f_j \|_{p }  \| \Box_{j,k-k_1-\ell}  g_i  \|_{\infty}
   .  \label{fgalgebra5}
\end{align}
By Young's inequality, it follows from \eqref{fgalgebra3} and \eqref{fgalgebra5} that
\begin{align}
 \|f_j g_i\|_{M^{[j]}_{p,1}}  \lesssim  \|f_j  \|_{M^{[j]}_{p,1}} \| g_i\|_{M^{[j]}_{\infty,1}} .  \label{fgalgebra6}
\end{align}
In view of Lemma \ref{Algebra1} and Bernstein's inequality,
\begin{align}
    \| g_i\|_{M^{[j]}_{\infty,1}} \lesssim  \| g_i\|_{M^{[i]}_{\infty,1}} \lesssim   2^{i d/p} \| g_i\|_{M^{[i]}_{p,1}} .  \label{fgalgebra7}
\end{align}
Hence, we have
\begin{align}
 \|f_j g_i\|_{M^{[j]}_{p,1}}  \lesssim  2^{id/p} \|f_j  \|_{M^{[j]}_{p,1}}  \| g_i\|_{M^{[i]}_{p,1}}.  \label{fgalgebra8}
\end{align}
It follows from \eqref{fgalgebra8} that
\begin{align}
  \|F_j\|_{M^{[j]}_{p,q}} \lesssim   \sum_{i\leq j}  \|f_j  \|_{M^{[j]}_{p,1}}  2^{id/p} \| g_i\|_{M^{[i]}_{p,1}}.  \label{fgalgebra9}
\end{align}
Hence,
\begin{align}
\left(\sum_{j\leq 0} 2^{j\mu r} \|F_j\|^r_{M^{[j]}_{p,q}}\right)^{1/r}
& \lesssim \left(\sum_{j\leq 0} 2^{j\mu r}  \|f_j  \|^r_{M^{[j]}_{p,1}} \right)^{1/r}  \sum_{i\leq 0}  2^{id/p} \| g_i\|_{M^{[i]}_{p,1}} \nonumber\\
& \lesssim   \|f\|_{\mathscr{M}^{\mu}_{p,1,r}} \|g\|_{\mathscr{M}^{d/p- \varepsilon}_{p,1,r}}  \label{fgalgebra10}
\end{align}
holds for any $0< \varepsilon \ll 1$. By the embedding $\mathscr{M}^{\mu}_{p,1,r} \subset \mathscr{M}^{d/p- \varepsilon}_{p,1,r}$ ($\mu \leq d/p-\varepsilon$) and \eqref{fgalgebra10}, we immediately have the result. In the case $r=1$, by the first inequality of \eqref{fgalgebra10}, we can obtain the desired estimate.  $\hfill\Box$\\

Since the decomposition to the frequency space in $\tilde{\mathscr{M}}^{\mu}_{p,1,r} $ is much rougher than that in $ \mathscr{M}^{\mu}_{p,1,r}$, it is harder to be an algebra for $\tilde{\mathscr{M}}^{\mu}_{p,1,r} $ than $ \mathscr{M}^{\mu}_{p,1,r}$. Indeed, we have

\begin{prop} \label{Algebra3}
Let $1\leq p<\infty$. Assume that $\mu \in [d/p, \max(d/p,d/p')]$ for $r=1$; and $p>p'$, $d/p < \mu \leq d/p'$ for $r\in (1,\infty)$. Then $\tilde{\mathscr{M}}^{\mu}_{p,1, r}$ is a Banach algebra.
\end{prop}

{\it Sketch Proof.} The proof is basically similar to that of Proposition \ref{Algebra2}.  Let $f,g\in \tilde{\mathscr{M}}^{\mu}_{p,1, r}$ with $f=\sum_{j\geq 0} f_j$ and $g=\sum_{i \geq 0} g_i$.  Further, we can assume that
$$
\left(\sum_{j\geq 0} 2^{j\mu r} \|f_j\|^r_{M^{[j]}_{p,1}} \right)^{1/r} \leq 2 \|f\|_{\tilde{\mathscr{M}}^{\mu}_{p,1,r}}, \ \ \ \left(\sum_{i  \geq 0} 2^{i \mu r} \|g_i\|^r_{M^{[i]}_{p,1}} \right)^{1/r} \leq 2 \|g\|_{\tilde{\mathscr{M}}^{\mu}_{p,1,r}}.
$$
We have
\begin{align}
fg   = \sum_{j\geq 0} \sum_{i\leq j} f_j g_i + \sum_{i\geq 0} \sum_{j<i} f_j g_i
  :=  \sum_{j\geq 0}  F_j +   \sum_{i\geq 0} G_i .  \label{3fgdecom}
\end{align}
Then we can repeat the procedures as in \eqref{fgalgebra2}--\eqref{fgalgebra9} to obtain that
\begin{align}
  \|F_j\|_{M^{[j]}_{p,q}} \lesssim  \|f_j  \|_{M^{[j]}_{p,1}}  \sum_{0\leq i\leq j} 2^{id/p} \| g_i\|_{M^{[i]}_{p,1}}.  \label{3fgalgebra9}
\end{align}
 By \eqref{3fgalgebra9} we can get the result, as desired. $\hfill\Box$

\begin{prop} \label{Algebra4}
Let $1\leq p<\infty$.   Then $\dot{\mathscr{M}}^{d/p}_{p,1, 1}$ is a Banach algebra.
\end{prop}
{\it Sketch Proof.}   Let $f,g\in \dot{\mathscr{M}}^{d/p}_{p,1, 1}$ with $f=\sum_{j\in \mathbb{Z}} f_j$ and $g=\sum_{i \in \mathbb{Z} } g_i$.  Further, we can assume that
$$
\sum_{j \in \mathbb{Z} } 2^{j d/p } \|f_j\|_{M^{[j]}_{p,1}}  \leq 2 \|f\|_{\dot{\mathscr{M}}^{d/p}_{p,1,1}}, \ \ \ \sum_{i \in \mathbb{Z} } 2^{i d/p} \|g_i\|_{M^{[i]}_{p,1}}  \leq 2 \|g\|_{\dot{\mathscr{M}}^{d/p}_{p,1,1}}.
$$
We have
\begin{align}
fg   = \sum_{j \geq 0 } \sum_{i \geq 0} f_j g_i + \sum_{j \leq 0} \sum_{i \leq 0} f_j g_i + \sum_{j \geq 0 } \sum_{i < 0} f_j g_i + \sum_{j \leq 0 } \sum_{i > 0} f_j g_i =I+II+III+IV.  \label{4fgdecom}
\end{align}
Using the same ways as in the proofs of Propositions \ref{Algebra2} and  \ref{Algebra3}, we can obtain that
\begin{align}
\|I\|_{\dot{\mathscr{M}}^{d/p}_{p,1,1}} + \|II\|_{\dot{\mathscr{M}}^{d/p}_{p,1,1}}  \lesssim   \|f\|_{\dot{\mathscr{M}}^{d/p}_{p,1,1}} \|g\|_{\dot{\mathscr{M}}^{d/p}_{p,1,1}}  \label{4fgalgebra10}
\end{align}
The estimates of  $\|III\|_{\dot{\mathscr{M}}^{d/p}_{p,1,1}}$ and $\|IV\|_{\dot{\mathscr{M}}^{d/p}_{p,1,1}}$ are similar. It suffices to consider the estimates of $\|IV\|_{\dot{\mathscr{M}}^{d/p}_{p,1,1}}$.  Taking $H_i=0 $ for $i\leq 0$, and $H_i = g_i \sum_{j \leq 0 } f_j$  for $i>0$,  we have
\begin{align}
 IV =  \sum_{i \in \mathbb{Z}} H_i = \sum_{i >0 } H_i.  \label{4fgdecom11}
\end{align}
Using the same way as in \eqref{fgalgebra2}--\eqref{fgalgebra9}, one obtains that for $i>0$
\begin{align}
  \|H_i\|_{M^{[i]}_{p,1}} \lesssim  \|g_i  \|_{M^{[i]}_{p,1}}  \sum_{ j \leq 0}   \| f_j \|_{M^{[i]}_{\infty,1}}  \lesssim  \|g_i  \|_{M^{[i]}_{p,1}}  \sum_{ j \leq 0} 2^{jd/p}  \| f_j \|_{M^{[j]}_{p,1}}.  \label{4fgalgebra9}
\end{align}
It follows from \eqref{4fgalgebra9} that
\begin{align}
\|IV\|_{\dot{\mathscr{M}}^{d/p}_{p,1,1}} \leq  \sum_{i>0}  2^{jd/p} \|H_i\|_{M^{[i]}_{p,1}}  \lesssim  \sum_{i>0}  2^{id/p} \|g_i  \|_{M^{[i]}_{p,1}}  \sum_{ j \leq 0} 2^{jd/p}  \| f_j \|_{M^{[j]}_{p,1}}.  \label{4fgalgebra10}
\end{align}
We obtain the result, as desired. $\hfill\Box$

\section{Local wellposedness of NLS in $\mathscr{M}^{\mu}_{p,1,r}$} \label{LNLS}

The solution of Schr\"odinger equation has a good behaviour on Feichtinger spaces $M^s_{p,q}$ (cf. \cite{BGOR,Wa06,BeOh20,ChFaSu12,CorNik08,CorNik092,KaKoIt14,WaHud07,WaHu07,WaHuHaGu11}). The local well-posedness of NLS in $M^0_{p,1}$ were obtained in \cite{Wa06,BeOk09}. We can generalize those local solutions in the scaling limit spaces $\mathscr{M}^{\mu}_{p,1,r}$.
We consider the Cauchy problem for a NLS
\begin{align} \label{NLS}
& {\rm i}u_t +\Delta u - F(u) =0,  \quad u(0,x)=u_0(x),
\end{align}
where $F(u) =\lambda |u|^{2\kappa} u$, $\kappa \in \mathbb{N}$; or $F(u)= (e^{\lambda|u|^2}-1) u$, $\lambda \in \mathbb{R}$.
Recall that  NLS has
the following equivalent form:
\begin{align}
& u(t)= S(t)u_0 +{\rm i}\mathscr{A} F(u), \label{NLSI}
\end{align}
where we denote
\begin{align}
 S(t)= e^{{\rm i} t\Delta}, \quad \mathscr{A}= \int^t_0 S(t-\tau) \cdot d\tau. \label{NLS1}
\end{align}
We have a local well-posed result for NLS in  $\mathscr{M}^{\mu}_{p,1,r}$.

\begin{thm} \label{NLSlocal}
Let $1\leq p,r<\infty$. Assume that $\mu \in [0, d/p)$ for $r>1$; $\mu \in [0,d/p]$ for $r=1$. Then NLS \eqref{NLSI} is local well-posed in $\mathscr{M}^{\mu}_{p,1,r}$. More precisely, for any $u_0 \in \mathscr{M}^{\mu}_{p,1,r}$, there exists a $T_m >0$ such that \eqref{NLSI} has a unique solution $u\in C([0, T_m); \mathscr{M}^{\mu}_{p,1,r})$ with
$$
\lim\sup_{t\to T_m} \|u(t)\|_{\mathscr{M}^{\mu}_{p,1,r}} =\infty.
$$
\end{thm}

By the embedding \eqref{4embedding2}, we see that $M^0_{p,1} \subset \mathscr{M}^{d/p}_{p,1} \subset M^0_{\infty,1}$.  Moreover, we have $\mathscr{M}^{d/p}_{p,1} \nsubseteq  M^0_{\tilde{p},1}$ for any $1\leq \tilde{p}<\infty$. It follows that $\mathscr{M}^{d/p}_{p,1} \nsubseteq  L^{\tilde{p}}$ for all $1\leq \tilde{p}<\infty$. Let $k_j =(j,0,...,0) \in \mathbb{Z}^d$, $f\in \mathscr{S}$  with ${\rm supp} \ \widehat{\varphi} \subset [-1/8,1/8]^d $ and
\begin{align} \label{examp1}
f= \sum_{j\leq -10} \frac{1}{j \ln ^2 |j|} e^{{\rm i} x k_j} \varphi(2^j x).
\end{align}
Then we have $f\in \mathscr{M}^{d/p}_{p,1}\setminus  M^0_{\tilde{p},1}$.  Indeed, let us observe that
\begin{align} \label{examp2}
\widehat{f} (\xi)= \sum_{j\leq -10} \frac{2^{-jd}}{j \ln ^2 |j|}   \widehat{\varphi}(2^{-j} (\xi-k_j) )
\end{align}
and $\Box_{k_j} f = \mathscr{F}^{-1} \psi(\cdot - k_j) \widehat{f}$. Since $\psi(\xi) =1$ for $|\xi|_\infty \leq 1/4$ and ${\rm supp} \ \psi \subset [-3/4, 3/4]^d$, we have
\begin{align} \label{examp3}
\widehat{\Box_{k_j} f} (\xi) =   \frac{2^{-jd}}{j \ln ^2 |j|}   \widehat{\varphi}(2^{-j} (\xi-k_j) ).
\end{align}
So, we have
\begin{align} \label{examp4}
\|f\|_{M^0_{\tilde{p}, 1}}  > \|\Box_{k_j} f \|_{\tilde{p}}  =   \frac{2^{-jd/\tilde{p}} }{|j| \ln ^2 |j|}  \|\varphi\|_{\tilde{p}}.
\end{align}
It follows that $\|f\|_{M^0_{\tilde{p}, 1}} =\infty$.  On the other hand, let $f_j=0$ for $j=0,-1,...,-9$ and
$$
f_j =  \frac{1}{j \ln ^2 |j|}  e^{{\rm i} x k_j} \varphi(2^j x), \ \ j\leq -10.
$$
It follows that $f=\sum_{j\leq 0} f_j$ and
\begin{align} \label{examp4}
\|f\|_{\mathscr{M}^{d/p}_{p,1, 1}}  \leq \sum_{j\leq 10} 2^{jd/p} \|f_j\|_{M^{[j]}_{p,1}}.
\end{align}
Let us observe that
$$
\Box_{j, 2^{-j}k_j} f_j = f_j, \ \  \Box_{j, k} f_j = 0 \  \mbox {if}  \  k\neq  2^{-j}k_j .
$$
Hence,  we have
$$
\|f_j\|_{M^{[j]}_{p,1}}  \lesssim    \frac{2^{-jd/p} }{|j| \ln ^2 |j|}.
$$
Then we obtain that
\begin{align} \label{examp5}
\|f\|_{\mathscr{M}^{d/p}_{p, 1}}  \lesssim  \sum_{j\leq 10}    \frac{1}{|j| \ln ^2 |j|}  <\infty.
\end{align}

Recall that the critical space of NLS for $F(u)= \lambda |u|^{2\kappa} u$ is $H^{d/2-1/\kappa}$ or $L^{d\kappa}$.  By the above argument, we see that the initial data like \eqref{examp1} is a class of supercritical data  for the local well-posedness of NLS.

\subsection{Boundedness of $S(t)$ in $\mathscr{M}^{\mu}_{p,q,r}$ }

The boundedness of $S(t)$ in $\mathfrak{M}^{\varrho}_{p,q,r}$ was obtained by Benyi and Oh \cite{BeOh20}. Using their result, we can show that  $S(t)$ is also bound in $\mathscr{M}^{\mu}_{p,q,r}$ and the following is a different approach from \cite{BeOh20}. By Young's inequality, we have
\begin{align}
  \|\Box_{j,k} S(t) u_0\|_p   \leq    \|\mathscr{F}^{-1} (\psi_{j,k} e^{it |\xi|^2}  ) \|_1  \|  u_0\|_p.  \label{bdd1}
\end{align}
Using Bernstein's multiplier estimate, we have
\begin{align}
  \|\Box_{j,k} S(t) u_0\|_p    \lesssim  \langle 2^{2j} t\rangle^{d/2}  \|  u_0\|_p.  \label{bdd2}
\end{align}
By Plancherel's identity,
\begin{align}
  \|\Box_{j,k} S(t) u_0\|_2    \leq  \|  u_0\|_2.  \label{bdd3}
\end{align}
Taking $p=1, \infty$ in \eqref{bdd2} and making interpolations between  \eqref{bdd2} and \eqref{bdd3}, we have for all $p\in [1,\infty]$,
\begin{align}
  \|\Box_{j,k} S(t) u_0\|_p    \lesssim  \langle 2^{2j} t\rangle^{d|1/2-1/p|}  \|  u_0\|_p.  \label{bdd4}
\end{align}
It follows from \eqref{bdd4} and the almost orthogonality of $\Box_{j,k}$ that for any $j\leq 0$, $1\leq p,q \leq \infty$,
\begin{align}
  \|  S(t) u_0\|_{M^{[j]}_{p,q}}    \lesssim  \langle 2^{2j} t\rangle^{d|1/2-1/p|}  \|  u_0\|_{M^{[j]}_{p,q}}.  \label{bdd5}
\end{align}
By the definition of $\mathscr{M}^{\mu}_{p,q,r}$, from \eqref{bdd5} we immediately have

\begin{prop} \label{BddS(t)}
Let $1\leq p, q\leq \infty$, $1\leq r< \infty$ and $\mu$ satisfy \eqref{dualcondition}
$$
0\wedge d\left(\frac{1}{q} -\frac{1}{p}\right)  \leq \mu(p,q)  \leq d\left(\frac{1}{q} +\frac{1}{p}-1\right).
$$
Then we have $S(t): \mathscr{M}^{\mu}_{p,q,r} \to \mathscr{M}^{\mu}_{p,q,r}$ and
\begin{align}
  \|  S(t) u_0\|_{\mathscr{M}^{\mu}_{p,q,r}}    \lesssim  \langle  t\rangle^{d|1/2-1/p|}  \|  u_0\|_{\mathscr{M}^{\mu}_{p,q,r}} .  \label{bdd6}
\end{align}
\end{prop}

\subsection{Local Well-Posedness in $\mathscr{M}^{\mu}_{p,1,r}$}

We now prove Theorem \ref{NLSlocal} and  only prove the result for the case $F(u)= (e^{\lambda|u|^2}-1) u$.  Considering the mapping
\begin{align}
\mathcal{T}:  u(t) \to  S(t)u_0 +{\rm i}\mathscr{A} F(u), \label{nlslocal1}
\end{align}
we will show that it is a contraction mapping in the space
\begin{align}
 X= \{u \in C([0,T]; \mathscr{M}^{\mu}_{p,1,r}): \|u\|_{C([0,T]; \mathscr{M}^{\mu}_{p,1,r})} \leq  M  \} \label{nlslocal2}
\end{align}
for some $0<T \leq 1$ and $M=2C \|u_0\|_{\mathscr{M}^{\mu}_{p,1,r}}$.   By Proposition \ref{BddS(t)},
\begin{align}
\|\mathcal{T}   u(t) \|_{\mathscr{M}^{\mu}_{p,1,r}}  \leq  C \|u_0\|_{\mathscr{M}^{\mu}_{p,1,r}}    + T \max_{t\in [0,T]}\| F(u(t))\|_{\mathscr{M}^{\mu}_{p,1,r}}.  \label{nlslocal3}
\end{align}
Since $\mathscr{M}^{\mu}_{p,1,r}$ is a Banach algebra, one has that
\begin{align}
 \| F(u)\|_{\mathscr{M}^{\mu}_{p,1,r}} &  \leq \sum^\infty_{m=1}  \frac{\lambda^m}{m\!} \| |u|^{2m} u \|_{\mathscr{M}^{\mu}_{p,1,r}} \nonumber\\
 &  \leq \sum^\infty_{m=1}  \frac{\lambda^m}{m\!}  C^{2m}\|u\|^{2m+1}_{\mathscr{M}^{\mu}_{p,1,r}}  \nonumber\\
  &  =  \left( \exp(\lambda C \|u\|^2_{\mathscr{M}^{\mu}_{p,1,r}}) -1 \right) \|u\|_{\mathscr{M}^{\mu}_{p,1,r}}.
  \label{nlslocal4}
\end{align}
Hence, for any $u\in X$,
\begin{align}
\|\mathcal{T} u \|_{C([0,T];\mathscr{M}^{\mu}_{p,1,r} ) }  \leq  C \|u_0\|_{\mathscr{M}^{\mu}_{p,1,r}}    + T \left( \exp(\lambda C M^2) -1 \right) M.  \label{nlslocal5}
\end{align}
Noticing that
$$
||u|^{2m} u - |v|^{2m} v | \leq (2m+1)  (|u|^{2m}  + |v|^{2m})     | u -  v |,
$$
we can use the same way as in \eqref{nlslocal5} to show that
\begin{align}
\|\mathcal{T}u -  \mathcal{T} v \|_{C([0,T];\mathscr{M}^{\mu}_{p,1,r} ) }  \leq     T \left( \exp(\lambda C M^2) -1 \right) \| u -  v \|_{C([0,T];\mathscr{M}^{\mu}_{p,1,r} ) }.    \label{nlslocal6}
\end{align}
By choosing $0<T<1$ such that
$$
T \left( \exp(\lambda C M^2) -1 \right)  = 1/100,
$$
we see that $\mathcal{T}$ is a contraction mapping in the space $X$. So, $\mathcal{T}$ has a fixed point $u\in X$, which is a unique solution of NLS.  We can extend the solution by considering the mapping
\begin{align}
\mathcal{T}_1 u(t)= S(t-T)u (T) +{\rm i}  \int^t_T S(t-\tau) F(u(\tau)) d\tau.  \label{NLSI1}
\end{align}
The the solution can be extended from $[0,T]$ to $[T, T_1]$.  Repeating this procedure we obtain a maximal $T_m$ satisfying $u \in C([0, T_m); \mathscr{M}^{\mu}_{p,1,r})$ and $T_m$ is a blowup time. This finishes the proof of Theorem \ref{NLSlocal}.

\section{Global wellposedness for NLS in $\mathscr{M}^{\mu}_{2,1,r}$} \label{GNLS}

In this section we consider the global solution of NLS \eqref{NLS} in the space  $\mathscr{M}^{\mu}_{2,1,r}$. By  establishing the refined Strichartz estimates adapted to the wavelet basis in the frequency spaces, we can show the global well-posedness of NLS for sufficiently small data in  $\mathscr{M}^{\mu}_{2,1,r} \supset M^0_{2,1}$, where the global well-posedness of NLS in $M^0_{2,1}$ was obtained in \cite{WaHud07}.

\subsection{Strichartz estimates }

For convenience, we denote
\begin{align} \label{Strichartz1}
\frac{2}{\gamma(p)} = d\left(\frac{1}{2}- \frac{1}{p}\right).
\end{align}

\begin{prop} \label{STRI}
Let $1\leq q,r<\infty$, $2\leq p \leq \infty$. Then the following inequalities
\begin{align} \label{Strichartz2}
\|S(t)u_0\|_{M^{[j]}_{p,q}}  \lesssim  \min (|t|^{-2/\gamma(p)}, 2^{4j/\gamma(p)} )   \| u_0\|_{M^{[j]}_{p',q}}
\end{align}
uniformly holds for all $j\in \mathbb{Z}.$
\end{prop}
{\it Proof.}  \eqref{Strichartz2} is obtained in \cite{BeOh20}. The following is an alternate proof.  By the basic decay of $S(t)$,
\begin{align} \label{Strichartz3}
\|\Box_{j,k} S(t)u_0\|_{p}  \lesssim   |t|^{-2/\gamma(p)}    \| u_0\|_{p'}.
\end{align}
On the other hand, using Young's and H\"older's inequalities, one has that
\begin{align} \label{Strichartz4}
\|\Box_{j,k} S(t)u_0\|_{p}  \leq      \| \psi_{j,k} e^{\mathrm{i}t|\xi|^2} u_0\|_{p'} \leq      \| \psi_{j,k} \|_{p/(p-2)} \|\hat{u}_0\|_{p} \lesssim    2^{4j/\gamma(p)}    \| u_0\|_{p'}.
\end{align}
Combining the estimates in \eqref{Strichartz3} and \eqref{Strichartz4}, we obtain
 \begin{align} \label{Strichartz5}
\|\Box_{j,k} S(t)u_0\|_{p}  \lesssim  \min( |t|^{-2/\gamma(p)},  2^{4j/\gamma(p)})    \| u_0\|_{p'}.
\end{align}
Substituting $u_0$ by $\sum_{|\ell|_\infty \lesssim 1}\Box_{j,k+\ell}u_0 $ and taking the sequence $\ell^q$ norms in both sides of \eqref{Strichartz5}, we have the result, as desired. $\hfill\Box$

\begin{lem} \label{STRII}
Let $2\leq p \leq \infty$, $\gamma \geq 2\vee \gamma(p)$, $(\gamma,p)\neq (2,\infty)$ if $d=2$.  Then the following inequalities
 \begin{align} \label{StrII1}
\|\Box_{j,k} \mathscr{A}f\|_{L^{\gamma}_tL^p_x(\mathbb{R}^{d+1})}  \lesssim    2^{4j(1/\gamma(p)-1/\gamma)}   \| f\|_{L^{\gamma'}_tL^{p'}_x(\mathbb{R}^{d+1})} .
\end{align}
uniformly holds for all $j\in \mathbb{Z}, \ k\in \mathbb{Z}^d.$
\end{lem}
{\it Proof.} By \eqref{Strichartz5}, we see that
\begin{align} \label{StrII2}
\|\Box_{j,k} \mathscr{A}f (t)\|_{p}  \lesssim  \int^t_0   \min( |t-\tau|^{-2/\gamma(p)},  2^{4j/\gamma(p)})    \| f(\tau)\|_{p'} d\tau.
\end{align}
Hence,
\begin{align} \label{StrII3}
 \|\Box_{j,k} \mathscr{A}f\|_{L^{\gamma}_tL^p_x(\mathbb{R}^{d+1})}  \lesssim  & 2^{4j/\gamma(p)} \left\| \int^t_0 \mathbbm{1}_{(|t-\tau|\leq 2^{-2j})}   \| f(\tau)\|_{p'} d\tau \right\|_{L^{\gamma}_t} \nonumber\\
 & +   \left\| \int^t_0  |t-\tau|^{-2/\gamma(p)} \mathbbm{1}_{(|t-\tau|> 2^{-2j})} \| f(\tau)\|_{p'} d\tau \right\|_{L^{\gamma}_t} \nonumber\\
 := & I+II.
\end{align}
By Young's inequality, we have
\begin{align} \label{StrII4}
 I   & \lesssim  2^{4j/\gamma(p)} \| \mathbbm{1}_{(|\cdot|\leq 2^{-2j})} \|_{\gamma/2}   \| f\|_{L^{\gamma'}_tL^{p'}_x} \nonumber\\
  & \lesssim    2^{4j(1/\gamma(p)-1/\gamma)}   \| f\|_{L^{\gamma'}_tL^{p'}_x}.
\end{align}
By changes of variables $t\to 2^{-2j} t$ and $\tau \to 2^{-2j}\tau$, we can find that
\begin{align} \label{StrII5}
 II  \lesssim    2^{4j/\gamma(p)-2j/\gamma-2j}   \left\| \int^t_0  |t-\tau|^{-2/\gamma(p)} \mathbbm{1}_{(|t-\tau|> 1)} \| f(2^{-2j}\tau)\|_{p'} d\tau \right\|_{L^{\gamma}_t}.
\end{align}
We divide the estimate of $II$ into the following three cases.

{\it Case 1.} If $\gamma(p) <2$, we see that $\gamma \geq 2$. It follows that $\mathbbm{1}_{(|t|\geq 1)} |t|^{-2/\gamma(p)} \in L^{\gamma/2}$.  By \eqref{StrII5} and Young's inequality, we have
\begin{align} \label{StrII6}
 II  \lesssim    2^{4j(1/\gamma(p)-1/\gamma)}   \| f\|_{L^{\gamma'}_tL^{p'}_x}.
\end{align}

{\it Case 2.} If $\gamma(p) >2$, we see that $\gamma \geq \gamma(p)$. In the case $\gamma> \gamma(p)$, then  $\mathbbm{1}_{(|t|\geq 1)} |t|^{-2/\gamma(p)} \in L^{\gamma/2}$. By \eqref{StrII5} and Young's inequality we see that \eqref{StrII6} holds. If $\gamma = \gamma(p)$, we can use Hardy-Littlewood-Sobolev's inequality to show that \eqref{StrII6} is true.

{\it Case 3.} If $\gamma(p) =2$, we see that $\gamma \geq 2$. In the case $\gamma> 2$, then  $\mathbbm{1}_{(|t|\geq 1)} |t|^{-1} \in L^{\gamma/2}$. By \eqref{StrII5} and Young's inequality we see that \eqref{StrII6} holds.

Summarizing the estimate of $I$ and those of $II$ in Cases 1--3, we see that we have shown \eqref{StrII1} holds except for the case $\gamma=\gamma(p)=2$.
If $\gamma=\gamma(p)=2$, using Keel and Tao's endpoint Strichartz estimates \cite{KeTa98}, we have
 \begin{align} \label{StrII7}
\| \mathscr{A}f\|_{L^{2}_tL^p_x }  \lesssim       \| f\|_{L^{2}_tL^{p'}_x },  \ \ p=\frac{2d}{d-2}, \ d\neq 2.
\end{align}
It follows that
 \begin{align} \label{StrII8}
\|\Box_{j,k} \mathscr{A}f\|_{L^{2}_tL^p_x(\mathbb{R}^{d+1})}  \lesssim    \| f\|_{L^{2}_tL^{p'}_x },  \ \ p=\frac{2d}{d-2}, \ d\neq 2.
\end{align}

\begin{lem} \label{STRIII}
Let $2\leq p \leq \infty$, $\gamma \geq 2\vee \gamma(p)$, $(\gamma,p)\neq (2,\infty)$ if $d=2$.  Then the following inequalities
 \begin{align} \label{StrIII1}
\|\Box_{j,k} S(t)u_0\|_{L^{\gamma}_tL^p_x(\mathbb{R}^{d+1})}  \lesssim    2^{2j(1/\gamma(p)-1/\gamma)}   \|u_0\|_{L^{2}(\mathbb{R}^{d})} .
\end{align}
uniformly holds for all $j\in \mathbb{Z}, \ k\in \mathbb{Z}^d.$
\end{lem}

{\it Proof.} By duality,
\begin{align} \label{StrIII2}
\left|\int_{\mathbb{R}} (\Box_{j,k}S(t)u_0, \psi(t)) dt\right| \leq \|u_0\|_2  \left\|\int_{\mathbb{R}} \Box_{j,k} S(-t)  \psi(t)  dt\right\|_2,
\end{align}
and by Lemma \ref{STRII},
\begin{align} \label{StrIII3}
  \left\|\int_{\mathbb{R}} \Box_{j,k} S(-t)  \psi(t)  dt\right\|^2_2 & \lesssim \|\psi\|_{L^{\gamma'}_tL^{p'}_x }  \left\|\int_{\mathbb{R}} \Box_{j,k} S(t-\tau)  \psi(\tau)  dt\right\| _{L^{\gamma}_tL^{p}_x } \nonumber\\
  & \lesssim  2^{4j(1/\gamma(p)-1/\gamma)}  \|\psi\|^2_{L^{\gamma'}_tL^{p'}_x }.
\end{align}
From \eqref{StrIII2} and \eqref{StrIII3} it follows that
\begin{align} \label{StrIII4}
\left|\int_{\mathbb{R}} (\Box_{j,k}S(t)u_0, \psi(t)) dt\right|  \lesssim  2^{2j(1/\gamma(p)-1/\gamma)}  \|u_0\|_2 \|\psi\|_{L^{\gamma'}_tL^{p'}_x }
\end{align}
By duality, \eqref{StrIII4} has implied the result of \eqref{StrIII1}. $\hfill\Box$

\begin{lem} \label{STRIV}
Let $2\leq p \leq \infty$, $\gamma \geq 2\vee \gamma(p)$, $(\gamma,p)\neq (2,\infty)$ if $d=2$.  Then the following inequalities
 \begin{align}
& \|\Box_{j,k} \mathscr{A}f\|_{L^{\infty}_tL^2_x(\mathbb{R}^{d+1})}  \lesssim    2^{2j(1/\gamma(p)-1/\gamma)}   \| f\|_{L^{\gamma'}_tL^{p'}_x(\mathbb{R}^{d+1})},  \label{StrIV1}\\
 & \|\Box_{j,k} \mathscr{A}f\|_{L^{\gamma}_tL^p_x(\mathbb{R}^{d+1})}  \lesssim    2^{2j(1/\gamma(p)-1/\gamma)}   \| f\|_{L^{1}_tL^{2}_x(\mathbb{R}^{d+1})} . \label{StrIV2}
\end{align}
uniformly hold for all $j\in \mathbb{Z}, \ k\in \mathbb{Z}^d.$
\end{lem}
{\it Proof.} The proof of \eqref{StrIV1} follows from \eqref{StrIII3} and \eqref{StrIV2} is the dual version of \eqref{StrIV1}. $\hfill\Box$

\begin{lem} \label{STRV}
Let $2\leq r, p \leq \infty$, $\gamma \geq 2\vee \gamma(p)$, $\beta\geq 2 \vee \gamma(r)$, and $(\gamma,p), (\beta,r) \neq (2,\infty)$ if $d=2$.  Then the following inequalities
\begin{align}
& \|\Box_{j,k} \mathscr{A}f\|_{L^{\beta}_tL^r_x (\mathbb{R}^{d+1})}  \lesssim    2^{2j(1/\gamma(p)-1/\gamma) +2j(1/\gamma(r)-1/\beta)}   \| f\|_{L^{\gamma'}_tL^{p'}_x(\mathbb{R}^{d+1})}   \label{StrV1}
\end{align}
uniformly holds for all $j\in \mathbb{Z}, \ k\in \mathbb{Z}^d.$
\end{lem}
{\it Proof.} In view of \eqref{StrII1} and \eqref{StrIV1}, using Berntein's inequality, we have for any $p\geq p_1$, $r_1\geq 2$,
\begin{align}
& \|\Box_{j,k} \mathscr{A}f\|_{L^{\gamma}_tL^{p_1}_x (\mathbb{R}^{d+1})}  \lesssim    2^{2j(1/\gamma(p)-1/\gamma) +2j(1/\gamma(p_1)-1/\gamma)}   \| f\|_{L^{\gamma'}_tL^{p'}_x(\mathbb{R}^{d+1})}   \label{StrV2} \\
& \|\Box_{j,k} \mathscr{A}f\|_{L^{\infty}_tL^{r_1}_x (\mathbb{R}^{d+1})}  \lesssim    2^{2j(1/\gamma(p)-1/\gamma) +2j /\gamma(r_1) }   \| f\|_{L^{\gamma'}_tL^{p'}_x(\mathbb{R}^{d+1})}   \label{StrV3}
\end{align}
We divide the proof into the following four cases.

{\it Case 1.} $\beta\in [\gamma, \infty]$ and $r\geq p\vee 2.$  Taking $p_1=r_1=r$ in \eqref{StrV2} and \eqref{StrV3}, $\theta = \gamma/\beta$. It follows that
$1/\beta = (1-\theta)/\infty + \theta/\gamma$. By H\"older's inequality,
\begin{align}
\|\Box_{j,k} \mathscr{A}f\|_{L^{\beta}_tL^r_x }  & \lesssim  \|\Box_{j,k} \mathscr{A}f\|^{1-\theta}_{L^{\infty}_tL^{r}_x }  \|\Box_{j,k} \mathscr{A}f\|^{\theta}_{L^{\gamma}_tL^{r}_x }  \nonumber \\
& \lesssim    2^{2j(1/\gamma(p)-1/\gamma) +2j(1/\gamma(r)-1/\beta)}   \| f\|_{L^{\gamma'}_tL^{p'}_x}.    \nonumber
\end{align}

{\it Case 2.} $\beta\in [\gamma, \infty]$ and $2\leq r \leq  p.$  Let $\theta = \gamma/\beta$, then $1/\beta = (1-\theta)/\infty + \theta/\gamma$.  One can choose suitable $p_1, r_1$ with $2\leq r_1 \leq r \leq p \leq p_1$ and $1/r= (1-\theta)/r_1 + \theta/p_1$ in \eqref{StrV2} and \eqref{StrV3}.
By H\"older's inequality,  and from \eqref{StrV2} and \eqref{StrV3} it follows that
\begin{align}
\|\Box_{j,k} \mathscr{A}f\|_{L^{\beta}_tL^r_x }  & \lesssim  \|\Box_{j,k} \mathscr{A}f\|^{1-\theta}_{L^{\infty}_tL^{r_1}_x }  \|\Box_{j,k} \mathscr{A}f\|^{\theta}_{L^{\gamma}_tL^{p_1}_x }  \nonumber \\
& \lesssim    2^{2j(1/\gamma(p)-1/\gamma) +2j(1/\gamma(r)-1/\beta)}   \| f\|_{L^{\gamma'}_tL^{p'}_x}.    \nonumber
\end{align}

{\it Case 3.} $\beta <\gamma$ and $ p  \geq 2\vee r.$ In view of \eqref{StrII1} and \eqref{StrIV2}, and the almost orthogonality of $\Box_{j,k}$, we have for any $p'_1\leq r'$, $r'_1 \leq 2$,
\begin{align}
 &\|\Box_{j,k} \mathscr{A}f\|_{L^{\beta}_tL^r_x }
  \lesssim    2^{2j(1/\gamma(r)-1/\beta) +2j(1/\gamma(p_1)-1/\beta)}   \| f\|_{L^{\beta'}_tL^{p'_1}_x},   \label{StrV4}\\
   &\|\Box_{j,k} \mathscr{A}f\|_{L^{\beta}_tL^r_x }
  \lesssim    2^{2j(1/\gamma(r)-1/\beta) +2j /\gamma(r_1) }   \| f\|_{L^{1}_tL^{r'_1}_x}.   \label{StrV5}
\end{align}
Since $p'\leq 2\wedge r'$, we see that \eqref{StrV4} and \eqref{StrV5} hold for $p'_1=r'_1=p'$.
Let $\theta = \beta/\gamma$, then $1/\gamma' = (1-\theta)  + \theta/\beta'$.
Interpolating \eqref{StrV4} and \eqref{StrV5} for $p'_1=r'_1=p'$, we immediately have \eqref{StrV1}.

{\it Case 4.} $\beta <\gamma$ and $ 2\leq p  \leq  r.$
It follows that $2 \geq p'\geq   r'$,
Let $\theta = \beta/\gamma$, then
\begin{align} \label{StrV6}
1/\gamma' = (1-\theta)  + \theta/\beta'.
\end{align}
Choosing suitable $r'_1$ and $p'_1$ such that $2\geq r'_1 \geq p' \geq r' \geq p'_1$ and
\begin{align} \label{StrV7}
1/p' = (1-\theta)/r'_1  + \theta/p'_1.
\end{align}
By  \eqref{StrV6} and \eqref{StrV7} we see that $(L^{1}_tL^{r'_1}_x, \ L^{\beta'}_tL^{p'_1}_x)_\theta = L^{\gamma'}_tL^{p'}_x$.
Making interpolation between \eqref{StrV4} and \eqref{StrV5}  , we obtain \eqref{StrV1}. $\hfill\Box$\\

For convenience, we denote
\begin{align}
\|f\|_{W^{[j]}_{\gamma,p,q}} = \left(\sum_{k\in \mathbb{Z}^d} \|\Box_{j,k} f\|^q_{L^\gamma_t L^p_x (\mathbb{R}^{d+1})} \right)^{1/q},  \label{Waveletmod}
\end{align}
and by $\mathscr{L}^{\gamma}(  \mathscr{M}^{\mu}_{p,q,r}) : = \mathscr{L}^{\gamma}( \mathbb{R},  \mathscr{M}^{\mu}_{p,q,r}) $ the space
\begin{align}
 \left\{f\in \mathscr{S}'(\mathbb{R}^{d+1}):  f=\sum_{j\leq 0} f_j, \   \left(\sum_{j\leq 0} 2^{j\mu r} \|f\|^r_{W^{[j]}_{\gamma,p,q}}\right)^{1/r}<\infty \right\}.    \nonumber
\end{align}
The norm on $\mathscr{L}^{\gamma}(  \mathscr{M}^{\mu}_{p,q,r})$ is defined by
\begin{align}
\|f\|_{\mathscr{L}^{\gamma}(  \mathscr{M}^{\mu}_{p,q,r})} =\inf \left(\sum_{j\leq 0} 2^{j\mu r} \|f\|^r_{W^{[j]}_{\gamma,p,q}}\right)^{1/r},
\end{align}
where the infimum is taken over all of the possible decompositions of $f=\sum_{j\leq 0} f_j \in \mathscr{L}^{\gamma}(  \mathscr{M}^{\mu}_{p,q,r})$. Similarly, we can define $\mathscr{L}^{\gamma}(  \dot{\mathscr{M}}^{\mu}_{p,q,r})$ by replacing the summation $\sum_{j\leq 0}$ with $\sum_{j\in \mathbb{Z}}$.

For our purpose we denote by $\mathscr{L}^{\gamma}(  \mathscr{M}^{\nu}_{p,q,r}) \dot{\cap} \mathscr{L}^{\infty}(  \mathscr{M}^{\mu}_{2,1,r})$ the following space
\begin{align}
 \left\{f\in \mathscr{S}'(\mathbb{R}^{d+1}):  f=\sum_{j\leq 0} f_j, \   \left(\sum_{j\leq 0} 2^{j\nu r} \|f\|^r_{W^{[j]}_{\gamma,p,q}}\right)^{1/r} \bigvee \left(\sum_{j\leq 0} 2^{j\mu r} \|f\|^r_{W^{[j]}_{\infty,2,1}}\right)^{1/r} <\infty \right\}.    \nonumber
\end{align}
The norm on  $\mathscr{L}^{\gamma}(  \mathscr{M}^{\nu}_{p,q,r}) \dot{\cap} \mathscr{L}^{\infty}(  \mathscr{M}^{\mu}_{2,1,r})$  is defined by
\begin{align}
\|f\|_{\mathscr{L}^{\gamma}(  \mathscr{M}^{\nu}_{p,q,r}) \dot{\cap} \mathscr{L}^{\infty}(  \mathscr{M}^{\mu}_{2,1,r})} =\inf \left\{ \left(\sum_{j\leq 0} 2^{j\nu r} \|f\|^r_{W^{[j]}_{\gamma,p,q}}\right)^{1/r}\bigvee \left(\sum_{j\leq 0} 2^{j\mu r} \|f\|^r_{W^{[j]}_{\infty,2,1}}\right)^{1/r} \right\},
\end{align}
where the infimum is taken over all of the possible decompositions of $f=\sum_{j\leq 0} f_j \in   \mathscr{L}^{\gamma}(  \mathscr{M}^{\nu}_{p,q,r}) \dot{\cap} \mathscr{L}^{\infty}(  \mathscr{M}^{\mu}_{2,1,r})$.

\begin{lem} \label{STRW}
Assume that the following conditions are satisfied:
\begin{itemize}
\item[\rm (i)] $2\leq p, p_1 \leq \infty$, $1\leq q\leq \infty$, $1\leq r<\infty$;

\item[\rm (ii)] $2\vee \gamma (p) \leq \gamma \leq \infty$, $2\vee \gamma (p_1) \leq \gamma_1 \leq \infty$;

\item[\rm (iii)] $0\wedge d(1/q-1/p) \leq \mu \leq d(1/p+1/q-1)$, $\delta (p,\gamma):= 2/\gamma(p) -2/\gamma$;
\item[\rm (iv)] $(\gamma, p), (\gamma_1, p_1) \neq (2, \infty)$ for $d=2$.
\end{itemize}
Then we have
\begin{align}
& \|S(t) u_0\|_{\mathscr{L}^{\gamma}(\mathscr{M}^{\mu-\delta (p,\gamma)}_{p,q,r})} \lesssim \|u_0\|_{\mathscr{M}^{\mu}_{2,q,r}}, \label{Strw1} \\
& \|\mathscr{A} f\|_{\mathscr{L}^{\gamma_1}(\mathscr{M}^{\mu-\delta(p_1,\gamma_1)}_{p_1,q,r})} \lesssim \|f\|_{\mathscr{L}^{\gamma'}(\mathscr{M}^{\mu+ \delta (p,\gamma) }_{p',q,r}) }.  \label{Strw2}
\end{align}
Moreover, if in addition $0   \leq \mu \leq d(1/p+1/q-1)$, then the above results also hold by replacing  $\mathscr{M}^{a}_{b,c,r}$ with  $\dot{\mathscr{M}}^{a}_{b,c,r}$.
\end{lem}
{\it Proof.}  \eqref{Strw1} follows from Lemma \ref{STRIII} and  \eqref{Strw2}  is a straightforward consequence of Lemma \ref{STRV}. $\hfill\Box$

\begin{lem} \label{STRX}
Let $p,q,r, \delta, \mu$  satisfy conditions of Lemma \ref{STRW}. Then we have
\begin{align}
& \|S(t) u_0\|_{\mathscr{L}^{\infty}(\mathscr{M}^{\mu }_{2,1,r}) \dot{\cap } \mathscr{L}^{\gamma}(\mathscr{M}^{\mu-\delta(p,\gamma)}_{p,1,r})} \lesssim \|u_0\|_{\mathscr{M}^{\mu}_{2,1,r}}, \label{Strw3} \\
& \|\mathscr{A} f\|_{\mathscr{L}^{\infty}(\mathscr{M}^{\mu }_{2,1,r}) \dot{\cap} \mathscr{L}^{\gamma}(\mathscr{M}^{\mu-\delta(p,\gamma)}_{p,1,r})} \lesssim \|f\|_{\mathscr{L}^{1}(\mathscr{M}^{\mu}_{2,1,r}) }.  \label{Strw4}
\end{align}
Moreover, if in addition $0   \leq \mu \leq d(1/p+1/q-1)$, then the above results also hold by replacing  $\mathscr{M}^{a}_{b,c,r}$ with  $\dot{\mathscr{M}}^{a}_{b,c,r}$.
\end{lem}
{\it Proof.} Taking $q=1$ in Lemmas \ref{STRIII} and \ref{STRV}, and noticing that $\delta(p,\gamma) =0$ for $(\gamma,p)= (\infty,2)$, we obtain \eqref{Strw3} and \eqref{Strw4}. $\hfill\Box$

\subsection{Nonlinear estiamtes}

\begin{lem} \label{NONEST2}
Let $d\geq 1$, $m\in \mathbb{N}$, $m\geq 2 \vee 4/d$,  $2/m =d(1/2-1/p)$, $1\leq r<\infty$. Assume that $0\leq \mu < d/2- 2/m$ for $r>1$; $0\leq \mu \leq d/2- 2/m$ for $r=1$.   Then we have
\begin{align}
\|u_1u_2...u_{m+1}\|_{\mathscr{L}^{1}(  \mathscr{M}^{\mu}_{2,1,r})  }  \lesssim  \prod^{m+1}_{\ell=1} \|u_\ell\|_{\mathscr{L}^{\infty}(  \mathscr{M}^{\mu}_{2,1,r}) \dot{\cap}  \mathscr{L}^{m}(  \mathscr{M}^{\mu}_{p,1,r})  }.  \label{Nonest2}
\end{align}
\end{lem}
{\it Proof.} Let $u_i \in \mathscr{L}^{m}(  \mathscr{M}^{\mu}_{p,1,r})  \dot{\cap}  \mathscr{L}^{\infty}(  \mathscr{M}^{\mu}_{2,1,r})$ and  $u_i$ is decomposed by
\begin{align}
 u_i =  \sum_{j_i \leq 0}   u_{i,j_i}, \ \ i=1,...,m+1.  \label{Nonest3}
\end{align}
We have
\begin{align}
u_1u_2...u_{m+1} =  \sum_{j_1,...,j_{m+1} \leq 0}   u_{1,j_1}u_{2,j_2}...u_{m+1, j_{m+1}}.  \label{Nonest4}
\end{align}
Let $\Lambda$ be the set of all permutations of $1,...,m+1$ and we denote by $\pi(1),...,\pi(m+1)$ the permutation of $1,...,m+1$.  We can rewrite \eqref{Nonest4} as
\begin{align}
u_1u_2...u_{m+1} = \sum_{\pi(1),...,\pi(m+1) \in \Lambda} \ \sum_{j_1\leq 0} \ \sum_{j_{m+1}\leq ...\leq j_1}   u_{\pi(1),j_1}u_{\pi(2),j_2}...u_{\pi(m+1), j_{m+1}}, \label{Nonest5}
\end{align}
where $\sum_{j_{m+1}\leq ...\leq j_1}:= \sum_{j_{m+1}=-\infty}^{ j_m}... \sum_{j_{2}=j_3}^{ j_1}$ denotes the summation on $j_2,...,j_{m+1}$.  By the symmetry, it suffices to estimate
\begin{align}
U =     \sum_{j_1\leq 0} \ \sum_{j_{m+1}\leq ...\leq j_1}   u_{1,j_1}u_{2,j_2}...u_{ m+1, j_{m+1}}: = \sum_{j_1\leq 0} U_{j_1}.   \label{Nonest6}
\end{align}
By the definition,
\begin{align}
\|U\|_{\mathscr{L}^{1}(  \mathscr{M}^{\mu}_{2,1,r})} \leq \left(\sum_{j_1 \leq 0} 2^{j_1\mu r} \|U_{j_1}\|^r_{W^{[j_1]}_{1,2,1}}\right)^{1/r}.  \label{Nonest7}
\end{align}
We estimate $\|U_{j_1}\|_{W^{[j_1]}_{1,2,1}}$. Using the almost orthogonality of $\Box_{j,k}$ and H\"older's inequality, we have
\begin{align}
&  \|\Box_{j_1,k} U_{j_1}\|_{L^1_t L^{2}_x} \nonumber\\
& \leq  \sum_{j_{m+1}\leq ...\leq j_1} \|\Box_{j_1,k}  (u_{1,j_1}u_{2,j_2}...u_{ m+1, j_{m+1}}) \|_{L^1_t L^{2}_x}  \nonumber\\
 & \leq  \sum_{j_{m+1}\leq ...\leq j_1} \sum_{k_1,..., k_{m+1} \in \mathbb{Z}^d } \|\Box_{j_1,k}  ( \Box_{j_1,k_1} u_{1,j_1} \, \Box_{j_1,k_2}u_{2,j_2}\, ... \, \Box_{j_1,k_{m+1}}u_{ m+1, j_{m+1}}) \|_{L^1_t L^{2}_x}  \nonumber\\
 & \leq  \sum_{j_{m+1}\leq ...\leq j_1} \sum_{k_1,..., k_{m+1} \in \mathbb{Z}^d } \| \Box_{j_1,k_1} u_{1,j_1} \|_{L^\infty_t L^{2}_x}  \prod^{m+1}_{i=2} \|\Box_{j_1,k_i}u_{i,j_i}\|_{L^m_t L^{\infty}_x} \mathbbm{1}_{(|k-k_1-...-k_{m+1}|_{\infty} \lesssim 1)} \label{Nonest8a}
\end{align}
By \eqref{Nonest8a} and  Young's inequality,
\begin{align}
\|U_{j_1}\|_{W^{[j_1]}_{1,2,1}}   & \leq  \sum_{j_{m+1}\leq ...\leq j_1}   \|  u_{1,j_1} \|_{W^{[j_1]}_{\infty, 2,1}}  \prod^{m+1}_{i=2} \| u_{i,j_i}\|_{W^{[j_1]}_{m, \infty, 1}} .
   \label{Nonest8}
\end{align}
It is easy to see that $W^{[j_2]}_{\gamma, \rho, 1}\subset W^{[j_1]}_{\gamma, \rho, 1}$ for $j_2\leq j_1 \leq 0$ and for $\rho_1 \geq \rho_2$,
\begin{align}
\|u\|_{W^{[j]}_{\gamma, \rho_1, q}} \lesssim 2^{j d(1/\rho_2 -1/\rho_1)}  \|u\|_{W^{[j]}_{\gamma, \rho_2, q}}.   \label{Nonest0}
\end{align}
Hence, it follows from \eqref{Nonest8} and \eqref{Nonest0} that
\begin{align}
 \| U_{j_1}\|_{W^{[j_1]}_{1,2, 1} }
& \leq  \sum_{j_{m+1}\leq ...\leq j_1}   \|  u_{1,j_1} \|_{W^{[j_1]}_{\infty,2,1}}  \prod^{m+1}_{i=2} \| u_{i,j_i}\|_{W^{[j_i]}_{m,\infty, 1}}  \nonumber\\
& \leq  \sum_{j_{m+1}\leq ...\leq j_1}   \|  u_{1,j_1} \|_{W^{[j_1]}_{\infty, 2,1}}  \prod^{m+1}_{i=2} 2^{j_i ( d/2-2/m)} \| u_{i,j_i}\|_{W^{[j_i]}_{m,p, 1}}
   \label{Nonest9}
\end{align}
Inserting the estimates in \eqref{Nonest9} into \eqref{Nonest7}, we obtain that
\begin{align}
\|U\|_{\mathscr{L}^{1}(  \mathscr{M}^{\mu}_{2,1,r})  }  \lesssim
 \left(\sum_{j_1\leq 0} 2^{j_1 \mu} \|  u_{1,j_1} \|^r_{W^{[j_1]}_{\infty, 2,1}} \right)^{1/r}  \prod^{m+1}_{i=2} \sum_{j_i\leq 0}  2^{j_i ( d/2-2/m)} \| u_{i,j_i}\|_{W^{[j_i]}_{m,p, 1}}.
\label{Nonest10}
\end{align}
Since $0\leq \mu < d/p = d/2-2/m$ for $r>1$, we have $\mathscr{L}^{m}(  \mathscr{M}^{d/p}_{p,1,1}) \supset \mathscr{L}^{m}(  \mathscr{M}^{\mu}_{p,1,r}) $. Hence, for $r>1$, it follows  from \eqref{Nonest10} that
\begin{align}
\|U\|_{\mathscr{L}^{1}(  \mathscr{M}^{\mu}_{2,1,r})  }  \lesssim     \prod_{  \ell=1,...,m+1} \|u_\ell\|_{\mathscr{L}^{ \infty}(  \mathscr{M}^{\mu}_{2,1,r}) \dot{\cap} \mathscr{L}^{m}(  \mathscr{M}^{\mu}_{p,1,r})  }.  \label{Nonest11}
\end{align}
The result for the case $r=1$ has been obtained in \eqref{Nonest10}.  $\hfill\Box$\\

In Lemma \ref{NONEST2}, we treated the power nonlinearity with condition $0\leq \mu \leq d/2-2/m$. If we consider an exponential nonlinear term $(e^{|u|^2} -1) u =\sum_{\kappa \geq 1} |u|^{2\kappa} u/\kappa !$, we use the following result.

\begin{lem} \label{NONESTE}
Let $d\geq 3$, $m\geq 2$, $p= 2d/(d-2)$,  $1\leq r<\infty$. Assume that $0\leq \mu < d/2- 1$ for $r>1$; $0\leq \mu \leq d/2- 1$ for $r=1$.   Then we have
\begin{align}
\|u_1u_2...u_{m+1}\|_{\mathscr{L}^{2}(  \mathscr{M}^{\mu}_{p',1,r})  }  \leq C^{m+1} \prod^{m+1}_{i=1}  \|u_i\|_{\mathscr{L}^{\infty}(  \mathscr{M}^{\mu}_{2,1,r})     \cap \mathscr{L}^{2}(  \mathscr{M}^{\mu}_{p,1,r})  }.  \label{ENonest2}
\end{align}
\end{lem}
{\it Sketch Proof.}
Let us connect the proof with \eqref{Nonest6}. By the definition,
\begin{align}
\|U\|_{\mathscr{L}^{2}(  \mathscr{M}^{\mu}_{p',1,r})} \leq \left(\sum_{j\leq 0} 2^{j_1\mu r} \|U_{j_1}\|^r_{W^{[j]}_{2,p',1}}\right)^{1/r}.  \label{2Nonest7}
\end{align}
We estimate $\|U_{j_1}\|_{W^{[j_1]}_{2,p',1}}$. Similar to \eqref{Nonest8a}
\begin{align}
&  \|\Box_{j_1,k} U_{j_1}\|_{L^2_t L^{p'}_x} \nonumber\\
& \leq  \sum_{j_{m+1}\leq ...\leq j_1} \sum_{k_1,..., k_{m+1} \in \mathbb{Z}^d } \| \Box_{j_1,k_1} u_{1,j_1} \|_{L^2_t L^{p}_x}  \prod^{m+1}_{i=2} \|\Box_{j_1,k_i}u_{i,j_i}\|_{L^\infty_t L^{dm/2}_x} \mathbbm{1}_{(|k-k_1-...-k_{m+1}|_{\infty} \lesssim 1)}.  \label{2Nonest8a}
\end{align}
By \eqref{2Nonest8a},
\begin{align}
\|U_{j_1}\|_{W^{[j_1]}_{2,p',1}} & \leq C m  \sum_{j_{m+1}\leq ...\leq j_1}   \|  u_{1,j_1} \|_{W^{[j_1]}_{2,p,1}}  \prod^{m+1}_{i=2} \| u_{i,j_i}\|_{W^{[j_1]}_{\infty, dm/2, 1}}.
   \label{2Nonest8}
\end{align}
Using \eqref{Nonest0} and similar to \eqref{Nonest9},
\begin{align}
\|U_{j_1}\|_{W^{[j]}_{2,p',1}} & \leq C^{m+1} \sum_{j_{m+1}\leq ...\leq j_1}   \|  u_{1,j_1} \|_{W^{[j_1]}_{2,p,1}}  \prod^{m+1}_{i=2} 2^{j_i (d/2-2/m)} \| u_{i,j_i}\|_{W^{[j_i]}_{\infty, 2, 1}}.
   \label{2Nonest9}
\end{align}
Then we can repeat the procedures as in the proof of Lemma \ref{NONEST2} to obtain the result, as desired. $\hfill\Box$

\subsection{Global well-posedness for NLS in $\mathscr{M}^{\mu}_{2,1,r}$}

We consider the following NLS:
\begin{align} \label{NLSp}
& {\rm i}u_t +\Delta u \pm |u|^{2\kappa} u=0,  \quad u(0,x)=u_0(x),
\end{align}
 \begin{thm} \label{NLSglobal}
Let $d\geq 1$, $\kappa\in \mathbb{N}$, $\kappa \geq 2/d$, $(\kappa, d) \neq (1,2)$,  $1\leq  r<\infty$. Assume that $\mu \in [0, d/2-1/\kappa)$ for $r>1$; $\mu \in [0,d/2-1/\kappa]$ for $r=1$. Then NLS \eqref{NLSp} is global well-posed for sufficiently small data in $\mathscr{M}^{\mu}_{2,1,r}$. More precisely, if $u_0 \in \mathscr{M}^{\mu}_{2,1,r}$  is sufficiently small,  then \eqref{NLSI} with $F(u) = \pm |u|^{2\kappa} u$ has a unique solution
$$
u\in C([0, \infty); \mathscr{M}^{\mu}_{2,1,r}) \cap  \left( \bigcap_{p\geq 2, \ \gamma\geq 2\vee \gamma(p)} \mathscr{L}^{\gamma} ( \mathscr{L}^{\mu-\delta(p,\gamma)}_{p,1,r}) \right).
$$
 \end{thm}

 By the embedding \eqref{4embedding1}, we see that $M^0_{2,1} \subset \mathscr{M}^{d/2-1/\kappa}_{2,1} \subset M^0_{d\kappa,1}$.  Moreover, we have $\mathscr{M}^{d/2-1/\kappa}_{2,1} \nsubseteq  M^0_{\tilde{p},1}$ for any $1\leq \tilde{p}<d\kappa $.  Let $k_j =(j,0,...,0) \in \mathbb{Z}^d$, $f\in \mathscr{S}$  with ${\rm supp} \ \widehat{\varphi} \subset [-1/8,1/8]^d $ and
\begin{align} \label{2examp1}
f= \sum_{j\leq -10} \frac{2^{j/\kappa}}{j \ln ^2 |j|} e^{{\rm i} x k_j} \varphi(2^j x).
\end{align}
Using the same way as in \eqref{examp1},  we can show that $f\in \mathscr{M}^{d/2-1/\kappa}_{2,1}\setminus  M^0_{\tilde{p},1}$ for $1\leq p < d\kappa$.   Recall that the global well-posedness of NLS in $M^0_{2,1}$ was established in \cite{WaHu07}, Theorem \ref{NLSglobal} is a generalization of the corresponding result of \cite{WaHu07}.

For an exponential nonlinearity $(e^{|u|^2} -1) u$,
\begin{align} \label{NLSe}
& {\rm i}u_t +\Delta u \pm (e^{|u|^2} -1) u =0,  \quad u(0,x)=u_0(x),
\end{align}
we have similar result:
 \begin{thm}
 \label{NLSglobal2}
Let $d\geq 3$,  $1\leq  r<\infty$. Assume that $\mu \in [0, d/2-1 )$ for $r>1$, $\mu \in [0,d/2-1 ]$ for $r=1$. Then NLS \eqref{NLSe} is global well-posed for sufficiently small data in $\mathscr{M}^{\mu}_{2,1,r}$. More precisely, for any sufficiently small $u_0 \in \mathscr{M}^{\mu}_{2,1,r}$,  \eqref{NLSI} with $F(u) = \pm (e^{|u|^2} -1) u$ has a unique solution
$$
u\in C([0, \infty); \mathscr{M}^{\mu}_{2,1,r}) \cap  \left( \bigcap_{p\geq 2, \gamma \geq 2\vee \gamma(p)} \mathscr{L}^{\gamma} ( \mathscr{L}^{\mu-\delta(p,\gamma)}_{p,1,r}) \right).
$$
\end{thm}

{\it Proof of Theoreom \ref{NLSglobal}.}  Let  $p$  satisfy $1/p = 1/2 -1/d\kappa $. Considering the mapping
\begin{align}
\mathcal{T}:  u(t) \to  S(t)u_0 +{\rm i}\mathscr{A} F(u), \label{nlsglobal1}
\end{align}
we will show that it is a contraction mapping in the space
\begin{align}
\mathscr{D}  = \{u \in \mathscr{L}^{2\kappa} ( \mathscr{M}^{\mu}_{p,1,r}) \cap \mathscr{L}^{\infty} ( \mathscr{M}^{\mu}_{2,1,r}) :  \|u\|_{\mathscr{L}^{2\kappa} ( \mathscr{M}^{\mu}_{p,1,r}) \cap \mathscr{L}^{\infty} ( \mathscr{M}^{\mu}_{2,1,r})} \leq M \} \label{nlsglobal2}
\end{align}
for  $M=2C \|u_0\|_{\mathscr{M}^{\mu}_{2,1,r}}$.  For simply, we write $X=\mathscr{L}^{2\kappa} ( \mathscr{M}^{\mu}_{p,1,r}) \cap \mathscr{L}^{\infty} ( \mathscr{M}^{\mu}_{2,1,r}) $.   By Corollary \ref{STRX},
\begin{align}
\|\mathcal{T}   u(t) \|_{X}  \leq  C \|u_0\|_{\mathscr{M}^{\mu}_{2,1,r}}  +C    \| F(u)\|_{\mathscr{L}^{1} ( \mathscr{M}^{\mu}_{2,1,r})}.  \label{nlsglobal3}
\end{align}
By Lemma \ref{NONEST2}
\begin{align}
  \| F(u)\|_{\mathscr{L}^{1} ( \mathscr{M}^{\mu}_{2,1,r})}  &  \leq    \|u\|^{2\kappa+1}_{X}.
  \label{nlsglobal4}
\end{align}
Hence, for any $u\in \mathscr{D}$,
\begin{align}
\|\mathcal{T} u \|_{X}  \leq  C \|u_0\|_{\mathscr{M}^{\mu}_{2,1,r}} +  C\|u\|^{2\kappa+1}_X .  \label{nlsglobal5}
\end{align}
Noticing that
$$
||u|^{2\kappa} u - |v|^{2\kappa} v | \leq (2\kappa+1)  (|u|^{2\kappa}  + |v|^{2\kappa})     | u -  v |,
$$
we can use the same way as in \eqref{nlsglobal5} to show that
\begin{align}
\|\mathcal{T}u -  \mathcal{T} v \|_{X}  \leq   C(\|u\|^{2\kappa}_X  + \|v\|^{2\kappa}_X )   \| u -  v \|_{X}.    \label{nlsglobal6}
\end{align}
By choosing  $M$ satisfying $CM^{2\kappa} \leq 1/100$,
we see that
$$
\|\mathcal{T} u \|_{X}  \leq M,  \ \ \ \|\mathcal{T}u -  \mathcal{T} v \|_{X}  \leq   \frac12  \| u -  v \|_{X}.
$$
Hence,
 $\mathcal{T}$ is a contraction mapping in the space $X$. So, $\mathcal{T}$ has a fixed point $u\in X$, which is a unique solution of NLS.  Using the embedding $\mathscr{L}^{\infty} ( \mathscr{M}^{\mu }_{2,1,r}) \subset C (\mathbb{R},  \mathscr{M}^{\mu }_{2,1,r})$, we see that $u\in C([0,\infty);  \mathscr{M}^{\mu }_{2,1,r})$. Again, in view of Corollary \ref{STRX}, we have $u\in \mathscr{L}^{\gamma} ( \mathscr{M}^{\mu-\delta(p,\gamma)}_{p,1,r})$ for any $\gamma \geq 2\kappa$. $\hfill\Box$\\

 {\it Sketch Proof of Theorem \ref{NLSglobal2}.} Following the proof of Theorem \ref{NLSglobal}, we substitute $\mathscr{D}$ in \eqref{nlsglobal2} by
 \begin{align}
\mathscr{E}  = \{u \in \mathscr{L}^{2} ( \mathscr{M}^{\mu}_{p,1,r}) \cap \mathscr{L}^{\infty} ( \mathscr{M}^{\mu}_{2,1,r}) :  \|u\|_{\mathscr{L}^{2} ( \mathscr{M}^{\mu}_{p,1,r}) \cap \mathscr{L}^{\infty} ( \mathscr{M}^{\mu}_{2,1,r})} \leq M \}, \label{2nlsglobal2}
\end{align}
where $p=2d/(d-2)$,  $M=2C \|u_0\|_{\mathscr{M}^{\mu}_{2,1,r}}$.  We write $\mathscr{Y}=\mathscr{L}^{2} ( \mathscr{M}^{\mu}_{p,1,r}) \cap \mathscr{L}^{\infty} ( \mathscr{M}^{\mu}_{2,1,r}) $.   By \eqref{Strw1} and \eqref{Strw2},
\begin{align}
\|\mathcal{T}   u(t) \|_{\mathscr{Y}}  \leq  C \|u_0\|_{\mathscr{M}^{\mu}_{2,1,r}}  +C    \| F(u)\|_{\mathscr{L}^{2} ( \mathscr{M}^{\mu}_{p',1,r})}.  \label{2nlsglobal3}
\end{align}
By Lemma \ref{NONESTE}, we have
\begin{align}
 \| F(u)\|_{\mathscr{Y}} &  \leq \sum^\infty_{m=1}  \frac{\lambda^m}{m\!} \| |u|^{2m} u \|_{\mathscr{Y}} \nonumber\\
 &  \leq \sum^\infty_{m=1}  \frac{\lambda^m}{m\!}  C^{2m}\|u\|^{2m+1}_{\mathscr{Y}}  \nonumber\\
  &  =  \left( \exp(\lambda C \|u\|^2_{\mathscr{Y}}) -1 \right) \|u\|_{\mathscr{Y}}.
  \label{2nlslocal4}
\end{align}
Then we can repeat the procedures as in the proof of Theorem \ref{NLSglobal} to obtain the result, as desired. $\hfill\Box$

\section{Appendix}

We show the following

\begin{prop}
Let $1\leq p,q<\infty$, $j\leq 0$. Then we have $(M^{[j]}_{p,q})^* = M^{[j]}_{p',q'}$ with equivalent norms and the equivalence is independent of $j\leq 0$.
\end{prop}
{\it Proof.} Similar to \eqref{dualarg1}, we have $M^{[j]}_{p',q'} \subset (M^{[j]}_{p,q})^*$ and
$$
\|g\|_{ (M^{[j]}_{p,q})^*} \leq C  \|g\|_{ M^{[j]}_{p',q'} },
$$
where $C$ is independent of $j\leq 0$ and $g \in M^{[j]}_{p',q'}$.  On the other hand, considering a mapping:
\begin{align}
\Gamma: M^{[j]}_{p,q} \ni f \mapsto \{\Box_{j,k} f\}_{k\in \mathbb{Z}^d} \in \ell^q(\mathbb{Z}^d, L^p(\mathbb{R}^d)),
\end{align}
we see that $\Gamma$ is an isometric mapping from $M^{[j]}_{p,q}$ onto
$$
\mathfrak{X}_0 := \left\{ \{\Box_{j,k} f\}_{k\in \mathbb{Z}^d}:  f\in  M^{[j]}_{p,q} \right\}
$$
and $g\in (M^{[j]}_{p,q})^*$ can be treated as a continuous linear functional in   $\mathfrak{X}_0 \subset \ell^q(\mathbb{Z}^d, L^p(\mathbb{R}^d))$. By Hahn-Banach's Theorem, $g$ can be extended onto   $ \ell^q(\mathbb{Z}^d, L^p(\mathbb{R}^d))$ and the extension is written as $\tilde{g}=\{g_k\}_{k\in \mathbb{Z}^d}$, whose norm will be preserved. Moreover, we have
$$
\langle \tilde{g}, \{\varphi_k\}_{k\in \mathbb{Z}^d} \rangle = \sum_{k\in \mathbb{Z}^d} \int_{ \mathbb{R}^d} \overline{\tilde{g}_k (x)} \varphi_k(x) dx.
$$
It follows that for any $f\in M^{[j]}_{p,q}$,
$$
\langle g, f \rangle = \langle \tilde{g}, \{\Box_{j,k} f \}_{k\in \mathbb{Z}^d} \rangle = \sum_{k\in \mathbb{Z}^d} \int_{ \mathbb{R}^d} \overline{\tilde{g}_k (x)} \Box_{j,k} f (x)  dx.
$$
Moreover, in view of Fubini's Theorem, we see that for any $f\in \mathscr{S}$,
$$
\langle g, f \rangle = \sum_{k\in \mathbb{Z}^d} \int_{ \mathbb{R}^d} \overline{\Box_{j,k}  \tilde{g}_k (x)} f (x)  dx.
$$
Hence, we have $g = \sum_{k\in \mathbb{Z}^d} \Box_{j,k}  \tilde{g}_k  $ and  by
$$
\left\|  \sum_{k\in \mathbb{Z}^d}  \Box_{j,k} \tilde{g}_k  \right\|_{ M^{[j]}_{p',q'}}  \lesssim  \left\|   \{ \tilde{g}_k \} \right\|_{\ell^q(\mathbb{Z}^d, L^p(\mathbb{R}^d)) } = \|g\|_{ (M^{[j]}_{p,q})^*}.
$$
\\

\noindent {\bf Acknowledgments.} The second named author is grateful to Professor Hans G. Feichtinger for his many enlightening discussions and for his kindly pointing out the references on Fofana spaces.  He is also grateful to Dr. Jie Chen for his carefully reading to the paper.    B.W. was supported in part by NSFC grant 1171024.

\end{document}